\documentclass[11pt,letterpaper]{article}
\usepackage{amsmath,amssymb,amsfonts,amsthm,euscript,changebar,multicol,amscd}
\usepackage{graphicx,wrapfig,subfigure,sidecap,epsfig,epstopdf}
\usepackage[table,dvipsnames]{xcolor}
\usepackage{hyperref}
\usepackage{url}
\usepackage[small,bf]{caption}
\usepackage{lastpage}
\usepackage[numbers,sort&compress]{natbib}
\usepackage{multirow}
\usepackage{cancel}
\usepackage{siunitx}
\usepackage[normalem]{ulem}
\usepackage{array}
\usepackage{varwidth}
\usepackage{lineno}

\hypersetup{
    colorlinks=true,       
    linkcolor=blue,       
    citecolor=blue,       
    filecolor=black,       
    urlcolor=black,         
    bookmarks=false,
    pdffitwindow=true,
    pdfpagelayout=SinglePage
}

\graphicspath{{figures/}}

\definecolor{bochkovcolor}{RGB}{255,191,0}

\definecolor{aslamcolor}{RGB}{150,150,150}

\newcommand{\mat}[1]{\mathbf{#1}}
\newcommand{\vect}[1]{\boldsymbol{#1}}

\newcommand{\Oof}[1]{\mathcal{O}\left( #1 \right)}
\newcommand{\ddt}[1]{\partial_t #1}
\newcommand{\ddn}[1]{\partial_{\vect{n}} #1}
\newcommand{\lap}{\nabla^2}
\newcommand{\hf}{\frac{1}{2}}

\usepackage[dvipsnames]{xcolor}

\newcommand{\reviewerOne}[1]{\textcolor{black}{#1}}
\newcommand{\reviewerTwo}[1]{\textcolor{black}{#1}}

\newcommand{\mylinelabel}[1]{}

\author{Daniil Bochkov, Frederic Gibou}
\title{PDE-Based Multidimensional Extrapolation of Scalar Fields over Interfaces with Kinks and High Curvatures}
\begin{document}
\maketitle
\begin{abstract}
We present a PDE-based approach for the multidimensional extrapolation of smooth scalar quantities across interfaces with kinks and regions of high curvature. \reviewerOne{Unlike the commonly used method of \cite{Aslam:04:A-partial-differenti}, in which normal derivatives are extrapolated, the proposed approach is based on the extrapolation and weighting of Cartesian derivatives. As a result, second- and third-order accurate extensions in the $L^\infty$ norm are obtained with linear and quadratic extrapolations, respectively, even in the presence of sharp geometric features.} The accuracy of the method is demonstrated on a number of examples in two and three spatial dimensions and compared to the approach of \cite{Aslam:04:A-partial-differenti}. \reviewerOne{The importance of accurate extrapolation near sharp geometric features is highlighted on an example of solving the diffusion equation on evolving domains.}
\end{abstract}

\section{Introduction}

Extrapolation procedures are ubiquitous in scientific computing and generally allow one to estimate a valid value of a quantity at points where data is not given; either in space or in time. In the context of level-set methods \cite{Osher;Sethian:88:Fronts-propagating-w}, extrapolation procedures in space have been frequently used since the advent of the ghost-fluid method \cite{Fedkiw;Aslam;Merriman;etal:99:A-Non-oscillatory-Eu}, where constant extrapolations were originally used. Generalized ghost-fluid methods were then designed, in part based on higher-order extrapolations for which Aslam introduced a partial differential equation (PDE) approach to perform linear and quadratic extrapolation \cite{Aslam:04:A-partial-differenti} and Gibou and Fedkiw introduced a cubic extrapolation in the same PDE framework \cite{Gibou;Fedkiw:05:A-fourth-order-accur}. It is natural in the level-set context to perform such extrapolations using PDE formulations for their solutions are based on Hamilton-Jacobi solvers that have been designed for other standard level-set equations, see e.g. \cite{Sussman;Smereka;Osher:94:A-level-set-approach}. A typical situation that needs extrapolation is that of an implicit treatment of a field in a free boundary problem. In this case, a valid value of the field at time $t^n$ needs to be known when assembling the right-hand side of the linear system of equations at time $t^{n+1}$. Since the interface at the new time step has swept grid points that are outside the domain at the previous time step, valid values of the field at time $t^n$ are needed in the domain at time $t^{n+1}$, which requires an extrapolation procedure.

Typical use of extrapolation methods can be found in a multitude of level-set applications including multiphase flow simulations \cite{Losasso;Gibou;Fedkiw:04:Simulating-Water-and, Gibou;Chen;Nguyen;etal:07:A-level-set-based-sh, Lepilliez;Popescu;Gibou;etal:16:On-two-phase-flow-so, Nguyen;Gibou;Fedkiw:02:A-Fully-Conservative, Langavant;Guittet;Theillard;etal:17:Level-set-simulation, Gibou;Min:12:Efficient-symmetric-, Robinson-Mosher;Shinar;Gretarsson;etal:08:Two-way-coupling-of-}, in the solution of Poisson-Boltzmann \cite{Mirzadeh;Theillard;Helgadottir;etal:12:An-Adaptive-Finite-D, Helgadottir;Gibou:11:A-Poisson-Boltzmann-} and Poisson-Nerntz-Planck equations \cite{Mirzadeh;Gibou:14:A-conservative-discr} for studying transport in ionic solutions, in heat and diffusion flow problems \cite{Gibou;Fedkiw:05:A-fourth-order-accur, Gibou;Fedkiw;Cheng;etal:02:A-second-order-accur, Bochkov;Gibou:19:Solving-the-Poisson-, Papac;Gibou;Ratsch:10:Efficient-symmetric-}, in the study of epitaxial growth and diblock-copolymer self-assembly used in the semi-conductor industry \cite{Papac;Helgadottir;Ratsch;etal:13:A-level-set-approach, ouaknin2016self}, in shape optimization \cite{Allaire;Jouve;Toader:03:Structural-optimizat, Theillard;Djodom;Vie;etal:13:A-second-order-sharp}, surface reconstruction of biomolecules \cite{Mirzadeh;Theillard;Helgadottir;etal:12:An-Adaptive-Finite-D, Egan:2018aa} and in Stefan-type problems \cite{Gibou;Fedkiw;Caflisch;etal:03:A-Level-Set-Approach, Chen;Min;Gibou:09:A-numerical-scheme-f, Mistani;Guittet;Bochkov;etal:17:The-Island-Dynamics-, Rycroft;Gibou:12:Simulations-of-a-str}. PDE-based extrapolation procedures have also been extended to adaptive Quad-/Oc-tree grids and parallel architectures \cite{Min;Gibou:07:A-second-order-accur, Mirzadeh;Guittet;Burstedde;etal:16:Parallel-level-set-m, Gibou;Min;Fedkiw:13:High-resolution-shar}. In addition, fast methods have been introduced for computationally efficient extrapolation procedures using the Fast Marching method, including parallel implementations \cite{Sethian2000aa, Sethian2001aa, Chacon;Vladimirsky:15:A-parallel-two-scale} or the Fast Sweeping method \cite{Zhao:04:A-Fast-Sweeping-Meth, ASLAM:2014aa}, including efficient parallel algorithms on adaptive grids \cite{Detrixhe;Gibou;Min:13:A-parallel-fast-swee, Detrixhe;Gibou:16:Hybrid-Massively-Par}. We also refer the interested reader to \cite{moroney2017extending} for another implicit approach to extrapolation based on solving the biharmonic equation.

However, those methods behave poorly in the case where the free boundary presents high-curvature features or kinks. \mylinelabel{rev1:kink_examples}\reviewerOne{Typical examples of such situations are multimaterial flows with triple junction points, motion of sharp-edged bodies in fluids, contact line dynamics in wetting phenomena, phase-change front propagation in the presence of confining walls, etc}. 
We introduce a method that solves that problem. We present the method in section \ref{sec:Numerical_Method} and numerical examples in sections \ref{sec:Numerical_Examples_2d} and \ref{sec:Numerical_Examples_3d} that illustrate its benefits and comment on its efficiency. \mylinelabel{rev1:example_intro}\reviewerOne{Section \ref{sec:Numerical_Diffusion} considers an example of solving a diffusion equation on evolving domains that demonstrates the importance of accurate extrapolation near sharp geometric features}. Section \ref{sec:Conclusion} draws some conclusions.

\section{Numerical Method} \label{sec:Numerical_Method}
\subsection{Level-set Representation}
The level set representation \cite{Osher;Sethian:88:Fronts-propagating-w} defines the interface of a domain by $\left\{\vect{x}: \phi(\vect{x}) = 0 \right\}$, its interior and exterior by $\phi(\vect{x}) < 0$ and $\phi(\vect{x}) > 0$, respectively, where $\phi(\vect{x})$ is a Lipschitz continuous function called \textit{the level-set function}. In this paper, the only geometrical quantity that is needed is the outward normal to the interface, $\vect{n}$, which can be computed as:
\begin{linenomath*}
\begin{eqnarray}
\vect{n} = \frac{\nabla \phi}{|\nabla \phi|}, \label{eq:Normal}
\end{eqnarray}
\end{linenomath*}
using central differencing for $\phi_x$ and $\phi_y$. In typical level-set simulations, the level-set function is reinitialized as a signed distance function \cite{Sussman;Smereka;Osher:94:A-level-set-approach}. We refer the interested reader to \cite{Sethian:96:Level-set-methods, Osher;Fedkiw:02:Level-Set-Methods-an} for a thorough presentation of the level-set method and \cite{Gibou:2018aa} for a recent review.
\subsection{\mylinelabel{rev1:name1}\reviewerOne{Normal-derivative based} multidimensional PDE extrapolation of \cite{Aslam:04:A-partial-differenti}}
High order extrapolations in the normal direction are traditionally performed in a series of steps, as proposed by Aslam in \cite{Aslam:04:A-partial-differenti} \mylinelabel{rev1:name2}\reviewerOne{and referred to in the present manuscript as the \textit{normal-derivative based partial differential equation (ND-PDE) extrapolation}}. For example, suppose that we seek to extrapolate a scalar field $q$ from the region where $\phi\leq 0$ to the region where $\phi>0$. In the case of a quadratic extrapolation, we first compute $q_{\vect{n}\vect{n}}=\nabla \left(\nabla  q \cdot \vect{n}\right)\cdot \vect{n}$ in the region $\phi \leq 0$ and
extrapolate it across the interface in a constant fashion, that is, such that its normal derivative is zero in the region $\phi > 0$, by solving the following partial differential equation:
\begin{linenomath*}
\begin{eqnarray}
\frac{\partial q_{\vect{n}\vect{n}}}{\partial \tau}+H(\phi) \left( \vect{n} \cdot \nabla q_{\vect{n}\vect{n}} \right)=0, \label{eq:Aslam1}
\end{eqnarray}
\end{linenomath*}
where $H$ is the Heaviside function. Then, the value of $q$ across the interface is found by solving the following two partial differential equations:
\begin{linenomath*}
\begin{align}
\label{eq:Aslam2}
\dfrac{\partial q_{\vect{n}}}{\partial \tau}  & +  H(\phi)\left(\vect{n} \cdot \nabla q_{\vect{n}} - q_{\vect{n}\vect{n}}\right) =  0,  \\
\dfrac{\partial q}{\partial \tau} & +  H(\phi)\left(\vect{n} \cdot \nabla q - q_{\vect{n}}\right)         =  0, \label{eq:Aslam3}
\end{align}
\end{linenomath*}
defining $q_{\vect{n}}$ in such a way that its normal derivative is equal to the previously extrapolated $q_{\vect{n}\vect{n}}$ and then defining $q$ in such a way that its normal derivative is equal to the previously extrapolated $q_{\vect{n}}$. These PDEs are solved in fictitious time $\tau$ for a few iterations (typically 15) since we only seek to extrapolate the values of $q$ in a narrow band of a few grid cells around the interface.

This extrapolation procedure produces accurate results in the case where the interface is smooth, but generates large error in the case where sharp geometric features occur, e.g. thin elongated shapes or interfaces with kinks as illustrated in sections \ref{sec:Numerical_Examples_2d} and \ref{sec:Numerical_Examples_3d}.

\subsection{\mylinelabel{rev1:name3}\reviewerOne{Weighted-Cartesian-derivative based} multidimensional PDE extrapolation}
Instead of calculating the normal derivatives in the negative region before extrapolating them, we instead compute the derivatives in the Cartesian directions, extrapolate them and then construct the normal derivatives. Specifically, consider the following quantities, that are computed in the negative level-set region:
\begin{linenomath*}
\begin{eqnarray*}
\vect{q}_{\nabla} =
\left(
\begin{array}{c}
 q_x  \\
 q_y   \\
 q_z  
\end{array}
\right) \textrm{ and the symmetric matrix } \mat{Q}_{\nabla\nabla} =
\left(
\begin{array}{ccc}
 q_{xx} & q_{xy} & q_{xz}  \\
 q_{xy} & q_{yy} & q_{yz}    \\
 q_{xz} & q_{zy} & q_{zz}   
\end{array}
\right).
\end{eqnarray*}
\end{linenomath*}
Similar to the method described in the previous section, we extrapolate the elements of $\mat{Q}_{\nabla\nabla}$ in a constant fashion:
\begin{linenomath*}
\begin{eqnarray}
\frac{\partial \mat{Q}_{\nabla\nabla}}{\partial \tau}+H(\phi)   \left(\vect{n} \cdot \nabla\mat{Q}_{\nabla\nabla}\right)=0,  \label{eq:Daniil1}
\end{eqnarray}
\end{linenomath*}
before successively solving the following equations:
\begin{linenomath*}
\begin{align}
\label{eq:Daniil2}
\dfrac{\partial \vect{q}_{\nabla}}{\partial \tau}  &+ H(\phi) \left( \vect{n} \cdot \left(\nabla \vect{q}_{\nabla} -  \mat{Q}_{\nabla\nabla} \right) \right) = 0, \\
\dfrac{\partial          q }{\partial \tau}  & +  H(\phi) \left( \vect{n} \cdot \left(\nabla q -  \vect{q}_{\nabla}\right)\right)  = 0.
\label{eq:Daniil3}
\end{align}
\end{linenomath*}
Note that now, the normal vector field $\vect{n}$ enters the equations merely as some sort of weighting factor. Thus, as long as field $q$ is sufficiently smooth this approach to multidimensional extrapolation is expected to produce accurate results even when the normal vector field $\vect{n}$ is not smooth (as is the case of domains with sharp features). \mylinelabel{rev1:name4}\reviewerOne{To distinguish the proposed approach from the one in \cite{Aslam:04:A-partial-differenti}, we refer to it as the \textit{weighted-Cartesian-derivative based partial differential equation (WCD-PDE)} extrapolation.}

\mylinelabel{rev2:higher-order}\reviewerTwo{\textbf{Remark:} It is possible to construct cubic and even higher-order extrapolations following this approach as well, however one needs to keep in mind the rapidly growing computational cost, because an $m$-th order method requires solving advection equations for tensor variables of order up to $m$ ($3\times 3 \times 3$ for cubic, $3\times 3\times 3 \times 3$ for quartic, etc).}

\subsection{Implementation details}

In this work we demonstrate the proposed method on uniform Cartesian grids and our implementation follows very closely the one from \cite{Aslam:04:A-partial-differenti} with just few differences. Consider a two dimensional computational grid with nodes defined as:
\begin{linenomath*}
\begin{align*}
  \vect{r}_{i,j} = 
  \begin{pmatrix}
    x_{\min} + (i-1) \Delta x \\
    y_{\min} + (j-1) \Delta y
  \end{pmatrix},
  \quad
  i \in [1;N_x],
  \,
  j \in [1;N_y],
  \quad
  \Delta x = \frac{x_{\max} - x_{\min}}{N_x-1},
  \quad
  \Delta y = \frac{y_{\max} - y_{\min}}{N_y-1},
\end{align*}
\end{linenomath*}
where $[x_{\min}; x_{\max}] \times [y_{\min}; y_{\max}]$ denotes the computational domain, $N_x$ and $N_y$ are number of grid nodes in the Cartesian directions. Standard second-order accurate central difference formulas are used for calculating the normal vector field $\vect{n}(\vect{r})$ (in the entire domain) and derivatives (first and second) of $q$ in the negative region. Normal derivatives of $q$ are computed as:
\begin{linenomath*}
\begin{align*}
  q_{\vect{n}} = \nabla q \cdot \vect{n} 
  \quad \textrm{and} \quad
  q_{\vect{n}\vect{n}} = \vect{n} \cdot \nabla\nabla q \cdot \vect{n} + \vect{n} \cdot \nabla\vect{n} \cdot \nabla q.
\end{align*}
\end{linenomath*}
Since the first and second order derivatives of $q$ are not well-defined at all grid points where $\phi < 0$ we replace the Heaviside function $\mathcal{H}(\phi)$ in equations \eqref{eq:Aslam2}, \eqref{eq:Daniil2} and in equations \eqref{eq:Aslam1}, \eqref{eq:Daniil1} with discrete fields $\mathcal{H}^{\phi, \nabla}$ and $\mathcal{H}^{\phi, \nabla\nabla}$, respectively, where:
\begin{linenomath*}
\begin{align*}
  \mathcal{H}^{\phi,\nabla}_{i,j} &=
  \begin{cases}
    0,\, \text{if } \phi \leq 0 \text{ at } \vect{r}_{i\pm 1, j},\,\vect{r}_{i,j\pm 1}, \\
    1,\, \textrm{otherwise},
  \end{cases}
  \\
  \mathcal{H}^{\phi,\nabla\nabla}_{i,j} &=
  \begin{cases}
    0,\,  \text{if } \phi \leq 0 \text{ at } \vect{r}_{i\pm 1, j},\,\vect{r}_{i,j\pm 1},\, \vect{r}_{i\pm 1, j \pm 1},\,\vect{r}_{i\pm 1,j\mp 1} \\
    1,\, \textrm{otherwise}.
  \end{cases}
\end{align*}
\end{linenomath*}
Applying an explicit first-order accurate in time discretization to equations \eqref{eq:Aslam1}-\eqref{eq:Daniil3} one obtains the following updating formulas:
\begin{linenomath*}
\begin{align}
\label{eq:Aslam_discete}
\begin{aligned}
\left[ q_{\vect{n}\vect{n}} \right]_{i,j}^{k+1} &= \left[ q_{\vect{n}\vect{n}} \right]_{i,j}^{k} - \Delta \tau \mathcal{H}^{\phi,\nabla\nabla}_{i,j} \left( \left[ \vect{n} \cdot \nabla q_{\vect{n}\vect{n}}\right]_{i,j}^k \right), \\
\left[ q_{\vect{n}} \right]_{i,j}^{k+1} &= \left[ q_{\vect{n}} \right]_{i,j}^{k} - \Delta \tau \mathcal{H}^{\phi,\nabla}_{i,j} \left( \left[ \vect{n} \cdot \nabla q_{\vect{n}}\right]_{i,j}^k - \left[q_{\vect{n}\vect{n}}\right]_{i,j}\right), \\
\left[ q \right]_{i,j}^{k+1} &= \left[ q \right]_{i,j}^{k} - \Delta \tau \mathcal{H}^{\phi}_{i,j} \left( \left[ \vect{n} \cdot \nabla q \right]_{i,j}^k - \left[q_{\vect{n}}\right]_{i,j}\right),
\end{aligned}
\end{align}
\end{linenomath*}
and
\begin{linenomath*}
\begin{align}
\label{eq:Daniil_discete}
\begin{aligned}
\left[ \mat{Q}_{\nabla\nabla} \right]_{i,j}^{k+1} &= \left[ \mat{Q}_{\nabla\nabla} \right]_{i,j}^{k} - \Delta \tau \mathcal{H}^{\phi,\nabla\nabla}_{i,j} \left( \left[ \vect{n} \cdot \nabla \mat{Q}_{\nabla\nabla} \right]_{i,j}^k \right), \\
\left[ \vect{q}_{\nabla} \right]_{i,j}^{k+1} &= \left[ \vect{q}_{\nabla} \right]_{i,j}^{k} - \Delta \tau \mathcal{H}^{\phi,\nabla}_{i,j} \left( \left[ \vect{n} \cdot \nabla \vect{q}_{\nabla} \right]_{i,j}^k - \left[\vect{n} \cdot \mat{Q}_{\nabla\nabla}\right]_{i,j}\right), \\
\left[ q \right]_{i,j}^{k+1} &= \left[ q \right]_{i,j}^{k} - \Delta \tau \mathcal{H}^{\phi}_{i,j} \left( \left[ \vect{n} \cdot \nabla q \right]_{i,j}^k - \left[\vect{n} \cdot \vect{q}_{\nabla}\right]_{i,j}\right),
\end{aligned}
\end{align}
\end{linenomath*}
\reviewerOne{
When extrapolating first- and second-order derivatives (i.e. $q_{\vect{n}}$, $q_{\vect{nn}}$, $\vect{q}_{\nabla}$ and $\mat{Q}_{\nabla\nabla}$), first-order spatial derivatives in the equations above are computed using first-order accurate upwind discretizations. For example, derivatives in the $x$-direction are approximated as:
\begin{linenomath*}
\begin{align}
\label{eq:1st-order-accurate-derivatives}
  \left[ n_x \partial_x f \right]_{i,j}^k = 
  \left\{
  \begin{aligned}
    \left[n_x\right]_{i,j} \frac{\left[f\right]^k_{i,j}-\left[f\right]^k_{i-1,j}}{\Delta x} + \Oof{\Delta x}, &\textrm{ if } \left[n_x\right]_{i,j} > 0, \\
    \left[n_x\right]_{i,j} \frac{\left[f\right]^k_{i+1,j}-\left[f\right]^k_{i,j}}{\Delta x} + \Oof{\Delta x}, &\textrm{ if } \left[n_x\right]_{i,j} < 0,
  \end{aligned}
  \right.
\end{align}
\end{linenomath*}
where $f$ is $q_{\vect{n}}$, $q_{\vect{nn}}$, $\vect{q}_{\nabla}$ or $\mat{Q}_{\nabla\nabla}$. This is sufficient to achieve second-order accuracy in the extended fields $q_{\vect{n}}$ and $\vect{q}_{\nabla}$. For extrapolation of the field $q$ itself (last equations in \eqref{eq:Aslam_discete} and \eqref{eq:Daniil_discete}), however, second-order accurate upwind discretizations are used. For example, derivatives in the $x$-direction are approximated as:
\begin{linenomath*}
\begin{align}
\label{eq:2nd-order-accurate-derivatives}
  \left[ n_x \partial_x q \right]_{i,j} = 
  \left\{
  \begin{aligned}
    \left[n_x\right]_{i,j} \left( \frac{\left[q\right]_{i,j}-\left[q\right]_{i-1,j}}{\Delta x} + \frac{\Delta x}{2} \textrm{minmod} \left( \left[q_{xx}\right]_{i,j}, \left[q_{xx}\right]_{i-1,j} \right) \right) + \Oof{\Delta x^2}, &\textrm{ if } \left[n_x\right]_{i,j} > 0, \\
    \left[n_x\right]_{i,j} \left( \frac{\left[q\right]_{i+1,j}-\left[q\right]_{i,j}}{\Delta x} - \frac{\Delta x}{2} \textrm{minmod} \left( \left[q_{xx}\right]_{i,j}, \left[q_{xx}\right]_{i+1,j} \right) \right) + \Oof{\Delta x^2}, &\textrm{ if } \left[n_x\right]_{i,j} < 0,
  \end{aligned}
  \right.
\end{align}
\end{linenomath*}
where
\begin{linenomath*}
\begin{align*}
  \textrm{minmod}(a,b) = 
  \begin{cases}
    0, \textrm{ if } ab \leq 0, \\
    a, \textrm{ if } |a| \leq |b|,\\
    b, \textrm{ if } |b| \leq |a|.
  \end{cases}
\end{align*}
\end{linenomath*}
Derivatives in the $y$-directions are approximated in a similar fashion. \mylinelabel{rev1:same}We note that approximation of derivatives 
is done in the same way for both, \reviewerOne{ND-PDE} and \reviewerOne{WCD-PDE}, extrapolation methods. The difference between the approaches lies in which quantities are extended over interfaces.}

\mylinelabel{rev1:numerics}\reviewerOne{Since in the new method the approximation of second-order derivatives in all Cartesian directions are already available during solving the PDE for $q$, the minmod corrections in \eqref{eq:2nd-order-accurate-derivatives}} can be computed only once during first iteration and reused in subsequent iterations, reducing the cost of each iteration by approximately 2 times. Specifically, the total count of arithmetic operations to compute $\left[ q \right]_{i,j}^{k+1}$ using the \mylinelabel{rev1:name5}\reviewerOne{ND-PDE} method is approximately 22 in two spatial dimensions and 32 in three spatial dimensions, while for the \mylinelabel{rev1:name6}\reviewerOne{WCD-PDE} method the total count is 10 and 14, correspondingly. Thus, if we denote as $T$ the cost of solving a single advection equation using first-order accurate approximations of derivatives, then the total cost of performing quadratic extrapolation using the \mylinelabel{rev1:name7}\reviewerOne{ND-PDE} method is approximately $(1+1+2)T=4T$ in two and three spatial dimensions, while the total cost of performing quadratic extrapolation using the \mylinelabel{rev1:name8}\reviewerOne{WCD-PDE} method is $(3+2+1)T=5T$ in two spatial dimensions and $(6+3+1)T=10T$ in three spatial dimensions.

\mylinelabel{rev1:stopping}\reviewerOne{Equations \eqref{eq:Aslam_discete} and \eqref{eq:Daniil_discete} are iterated in the fictitious time $\tau$ until steady-state. Time step $\Delta \tau$ is chosen based on consideration of satisfying the CFL condition as:
\begin{linenomath*}
\begin{align*}
  \Delta \tau = \frac{\min (\Delta x, \Delta y)}{2}
  \quad \textrm{and} \quad 
  \Delta \tau = \frac{\min (\Delta x, \Delta y, \Delta z)}{3}
\end{align*}
\end{linenomath*}
in two and three spatial dimensions, respectively. Iterations are terminated when the maximum difference between two successive steps $\max\limits_{i,j} \left| \left[ f \right]_{i,j}^{k+1} - \left[ f \right]_{i,j}^{k} \right|$ within the band of interest, that is, among all grid nodes within the distance of $2\sqrt{\Delta x^2 + \Delta y^2}$ (or $2\sqrt{\Delta x^2 + \Delta y^2  + \Delta z^2}$ in three spatial dimensions) around the domain boundary, is less than a specified tolerance $\epsilon_\textrm{tol} = 10^{-12}$.}

\textbf{Remark.} Since in the proposed approach there is no need to recalculate second derivatives and apply the nonlinear minmod operator at every iteration, it is possible to obtain the steady-state solution of the advection equations in an implicit fashion. This could be very beneficial in cases when a good guess for the extended field is available (for example, solutions from preceding time instants in time-dependent problems). Such an approach will be explored in future works.

\reviewerOne{
\mylinelabel{rev1:numerics2} \textbf{Remark.} In case of linear extrapolations, the first equations in \eqref{eq:Aslam_discete} and \eqref{eq:Daniil_discete} are not solved, in second equations $\left[ q_{\vect{nn}} \right]_{i,j}$ and $\left[ \mat{Q}_{\nabla\nabla} \right]_{i,j} $ are set to zero and first-order accurate formulas \eqref{eq:1st-order-accurate-derivatives} are used during the extrapolation of both field $q$ and its derivatives for more efficient computations.}

\subsubsection{Extension to adaptive Quad-/Oc-tree grids}

\begin{figure}[h]
\begin{center}
\includegraphics[height=.3\textwidth]{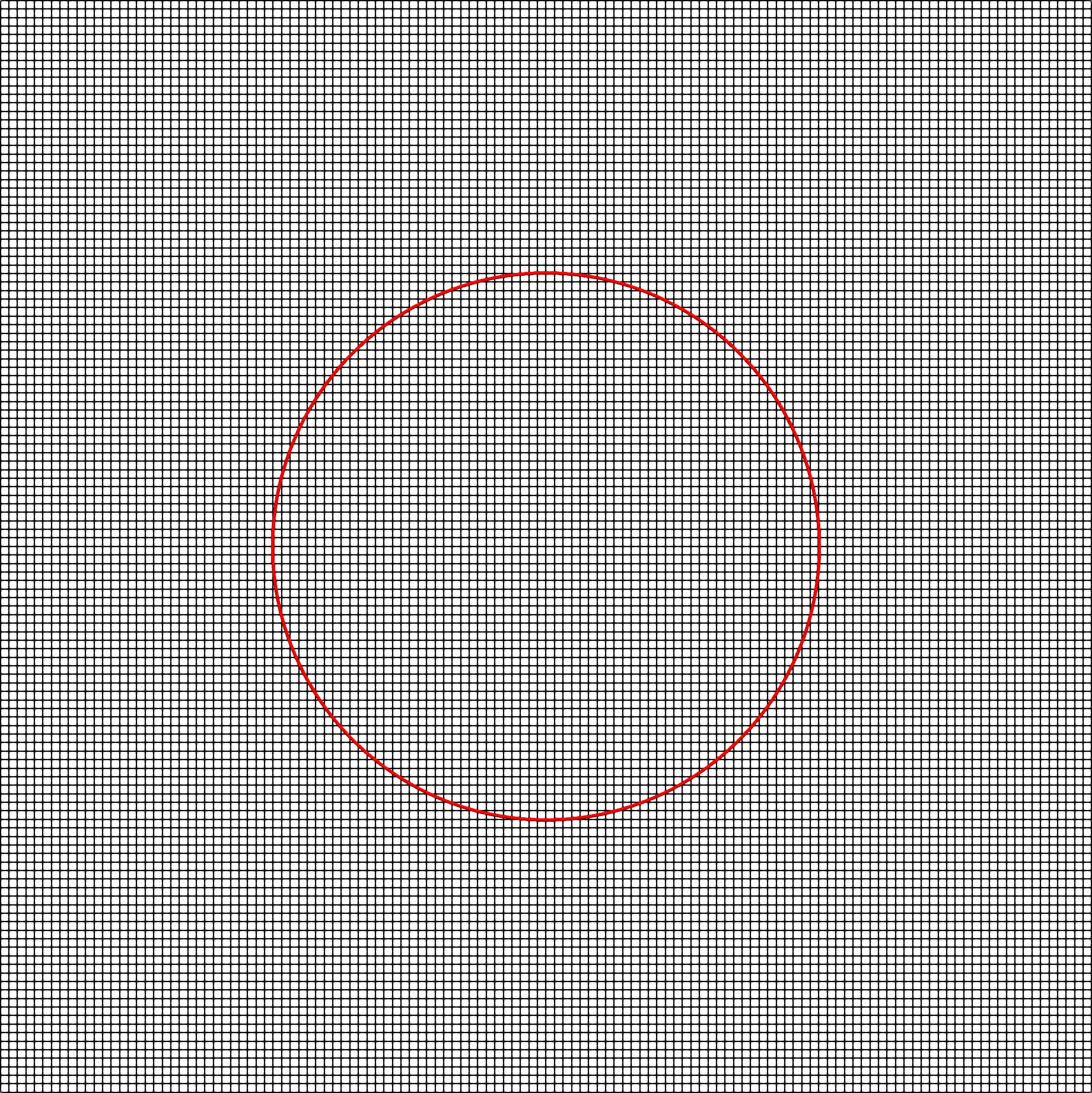}
$\quad$
\includegraphics[height=.3\textwidth]{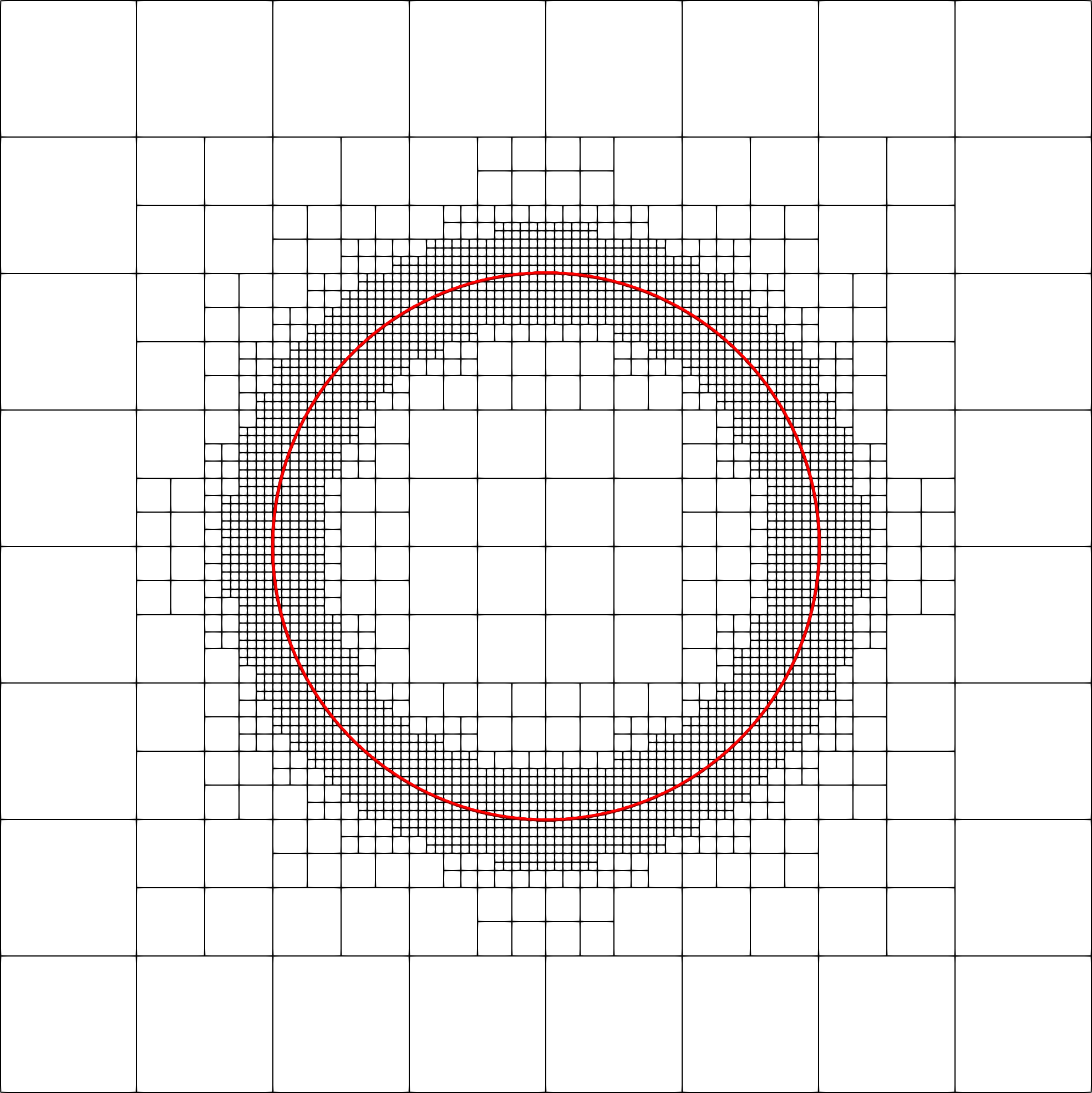}
\end{center}
\vspace{-.3cm}\vspace{-.3cm}\caption{\it Examples of uniform (left) and adaptive (right) Cartesian grids for a circular interface.} \label{fig::grids}
\end{figure}

The methodology introduced in this paper can be trivially extended to Quad-/Oc-tree data structures. Specifically, we sample data fields at nodes of Quad-/Oc-tree grids and use the second-order accurate discretizations of \cite{Min;Gibou:07:A-second-order-accur} for regions where grids are non-uniform. A band of uniform grids are usually imposed near the interface in practical free boundary applications (see Fig. \ref{fig::grids}). In this case the extrapolation within some neighborhood around the interface (where it is primarily required) is as accurate as for uniform grids, however, the extrapolation procedure is much faster on adaptive grids for their significant reduction in the total number of grid points.

\section{Numerical Results in Two Spatial Dimensions} \label{sec:Numerical_Examples_2d}

\begin{figure}[!h]
\begin{center}
\subfigure[Disk, $\Omega_0$]{{\includegraphics[height=.22\textwidth]{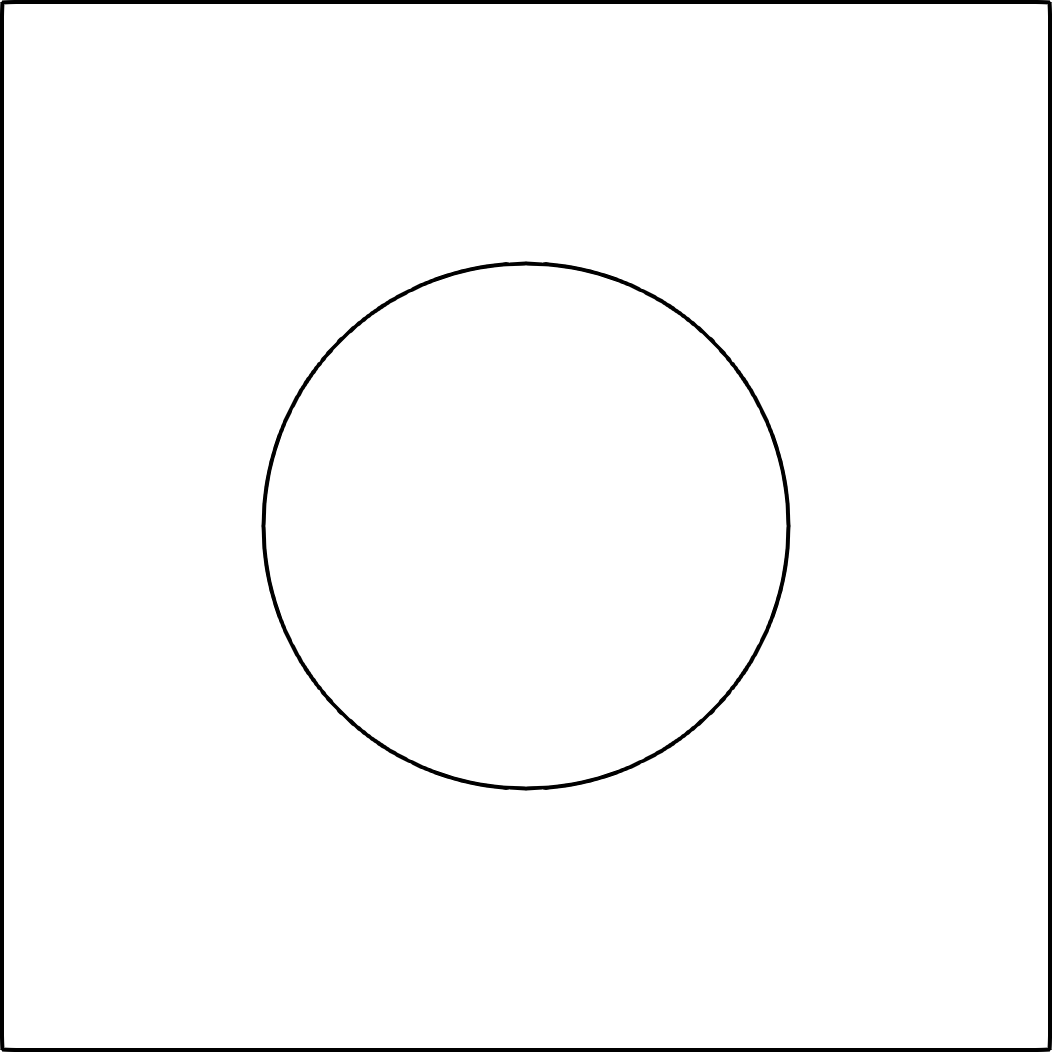}\label{fig:domain2_Aslam}}}  
\subfigure[Star, $\Omega_1$]{{\includegraphics[height=.22\textwidth]{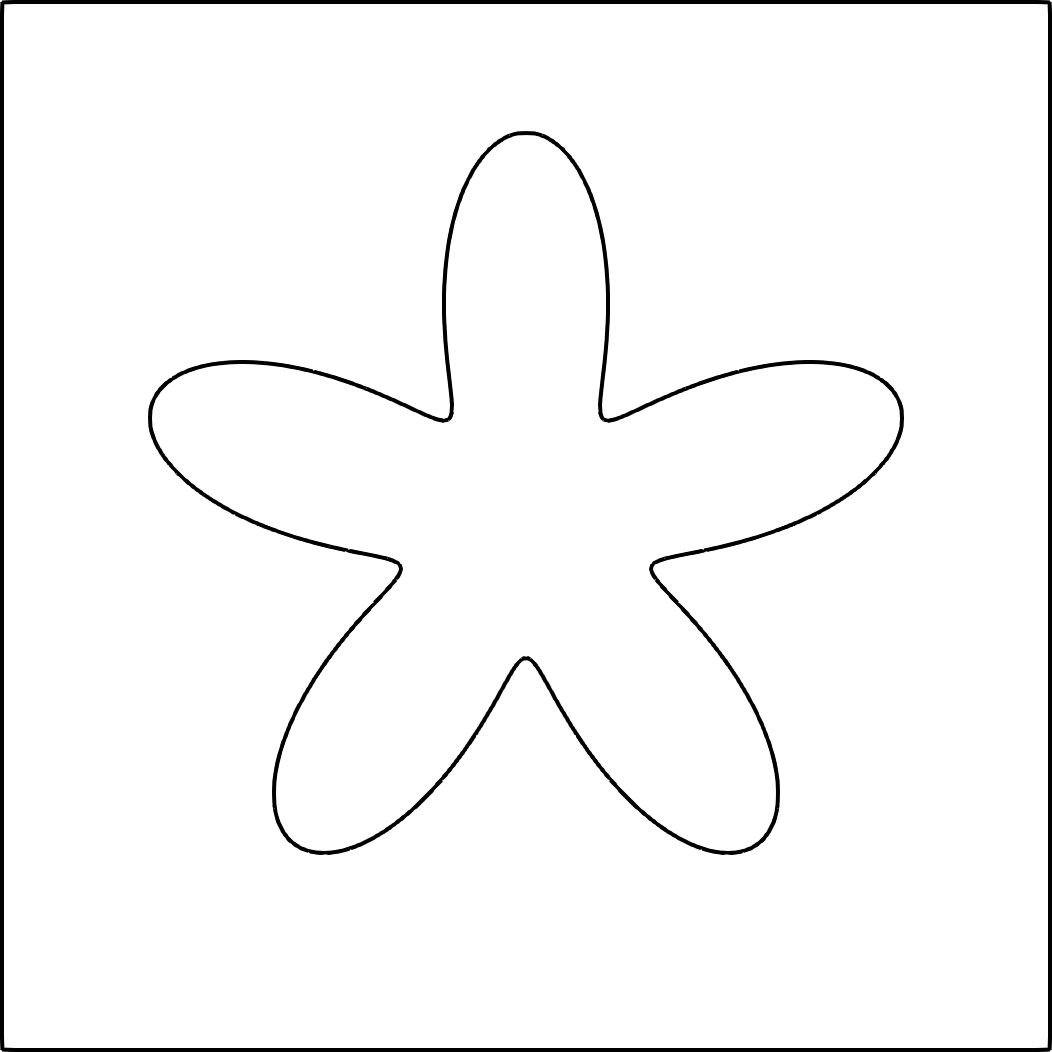}\label{fig:domain3_Aslam}}} 
\subfigure[Union, $\Omega_2$]{{\includegraphics[height=.22\textwidth]{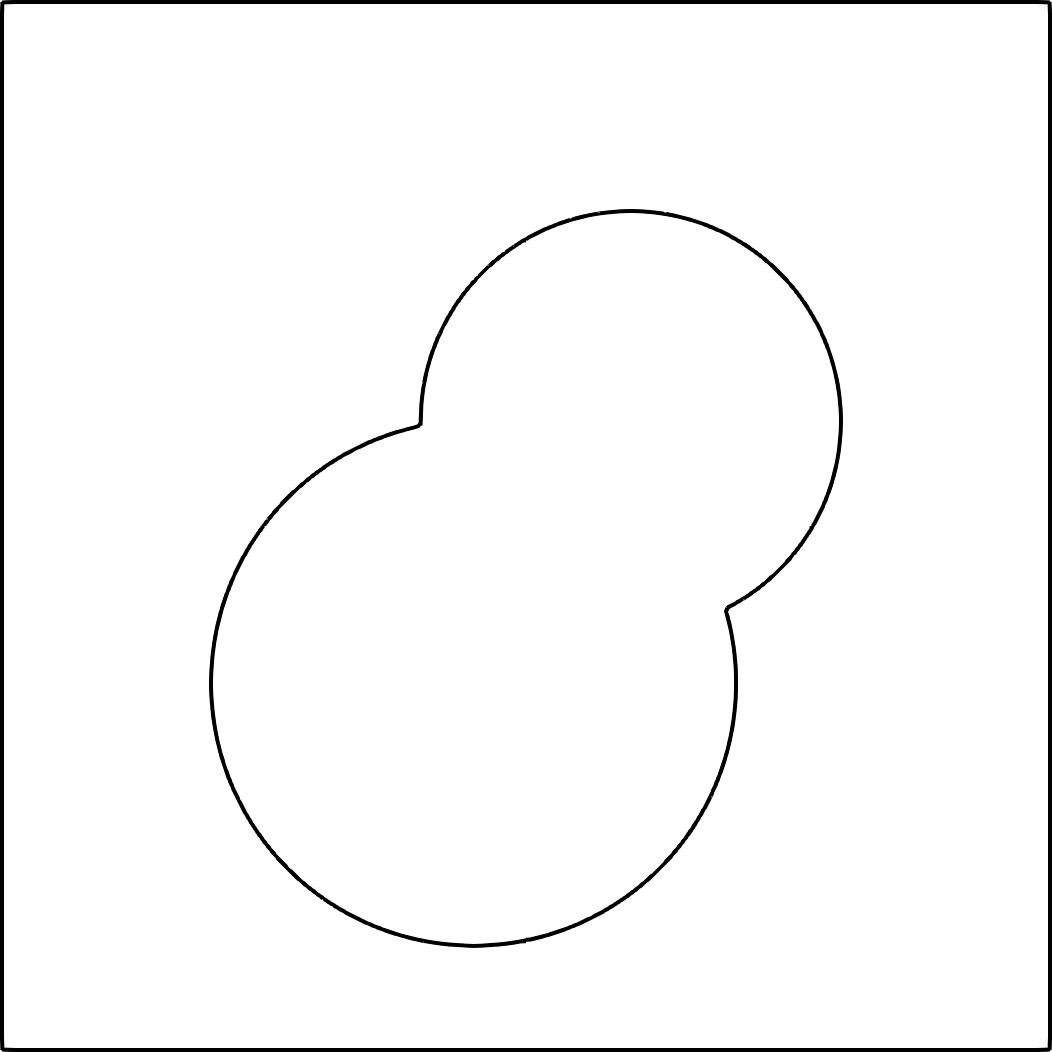}\label{fig:union}}} 
\subfigure[Intersection, $\Omega_3$]{{\includegraphics[height=.22\textwidth]{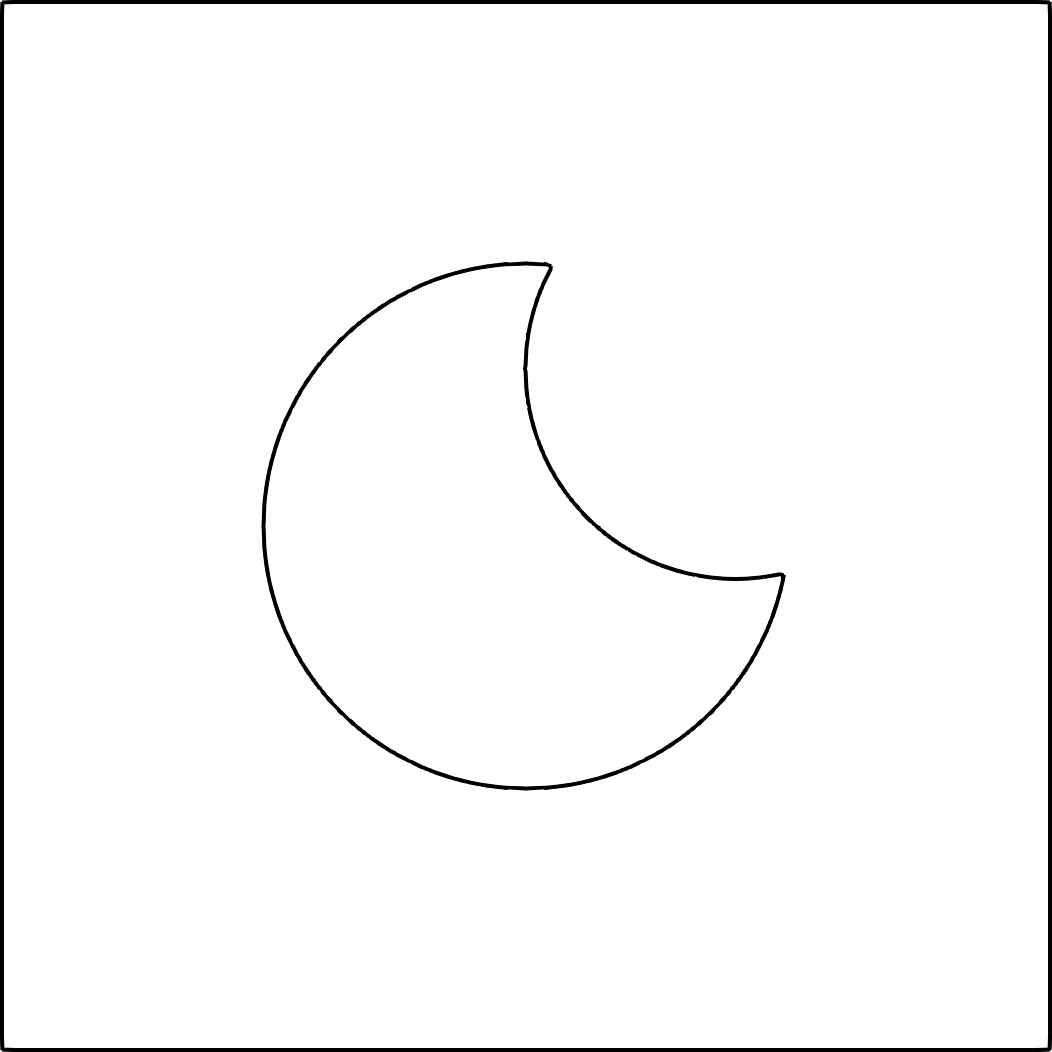}\label{fig:intersection}}} 
\end{center}
\vspace{-.3cm}\vspace{-.3cm}\caption{\it Computational domains considered in section \ref{sec:Numerical_Examples_2d}.} \label{fig::domains}
\end{figure}

We consider four physical domains: a disk, a star shape, a union of two disks and an intersection of two disks (see figure \ref{fig::domains}). The disk is a smooth interface for which the approach of \cite{Aslam:04:A-partial-differenti} performs well. The star-shape domain is an example where regions of high curvature are present (crest and trough of the wavy shape). The union/intersection of two disks are examples where kinks occur and illustrate the case of typical free boundary simulations where changes in topology occur. The definition of those domains are given by the level-set functions:
\begin{linenomath*}
\begin{align*}
\begin{aligned}
\phi_0(x, y) &= \sqrt{x^2 + y^2}-0.501, \\
\phi_1(x, y) &= \sqrt{x^2 + y^2}-0.501 - 0.25 \frac{y^5 + 5x^4y-10x^2y^3}{\left(x^2 + y^2\right)^{\frac{5}{2}}}, \\
\phi_2(x, y) &= \min\left(\sqrt{(x+.1)^2 + (y+.3)^2}-0.501,
                   \sqrt{(x-.2)^2 + (y-.2)^2}-0.401\right), \\
\phi_3(x, y) &= \max\left(\sqrt{x^2 + y^2}-0.501,
                   \sqrt{(x-.4)^2 + (y-.3)^2}-0.401\right),
\end{aligned}
\end{align*}
\end{linenomath*}

\begin{figure}[!h]
\begin{center}
\includegraphics[width=.99\textwidth]{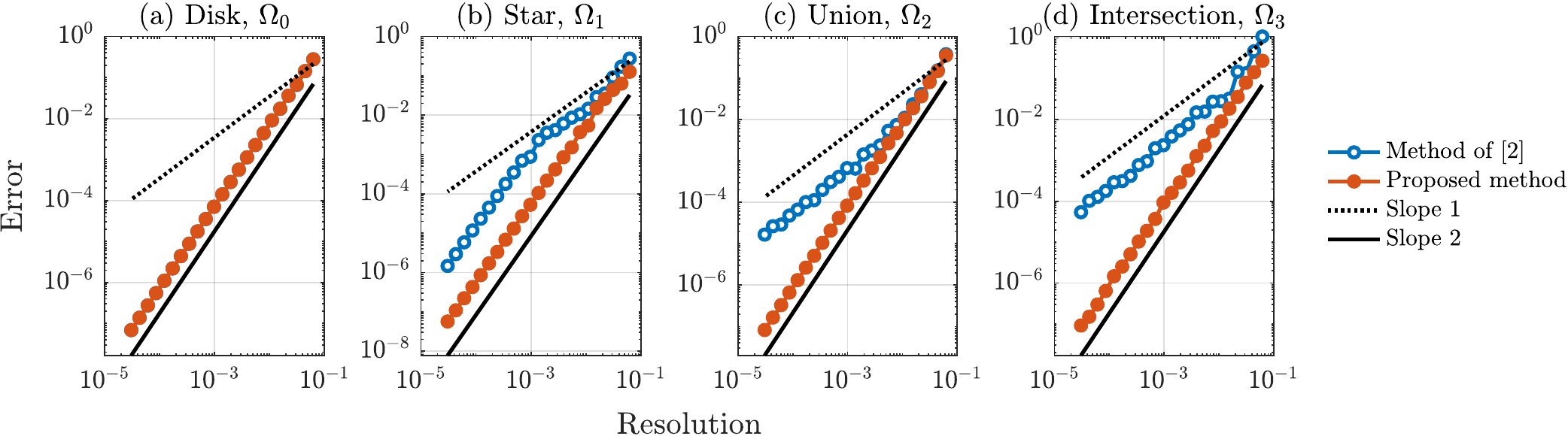}
\end{center}
\vspace{-.3cm}\vspace{-.3cm}\caption{\it Accuracy of the linear extrapolation \reviewerOne{(in the $L^\infty$ norm)} in two spatial dimensions measured in a narrow band of thickness $2\sqrt{\Delta x^2 + \Delta y^2}$ around an interface using \reviewerOne{the approach of \cite{Aslam:04:A-partial-differenti} and the proposed approach}.} \label{fig::results::2d_linear}
\end{figure}

\begin{figure}[!h]
\begin{center}
\includegraphics[width=.99\textwidth]{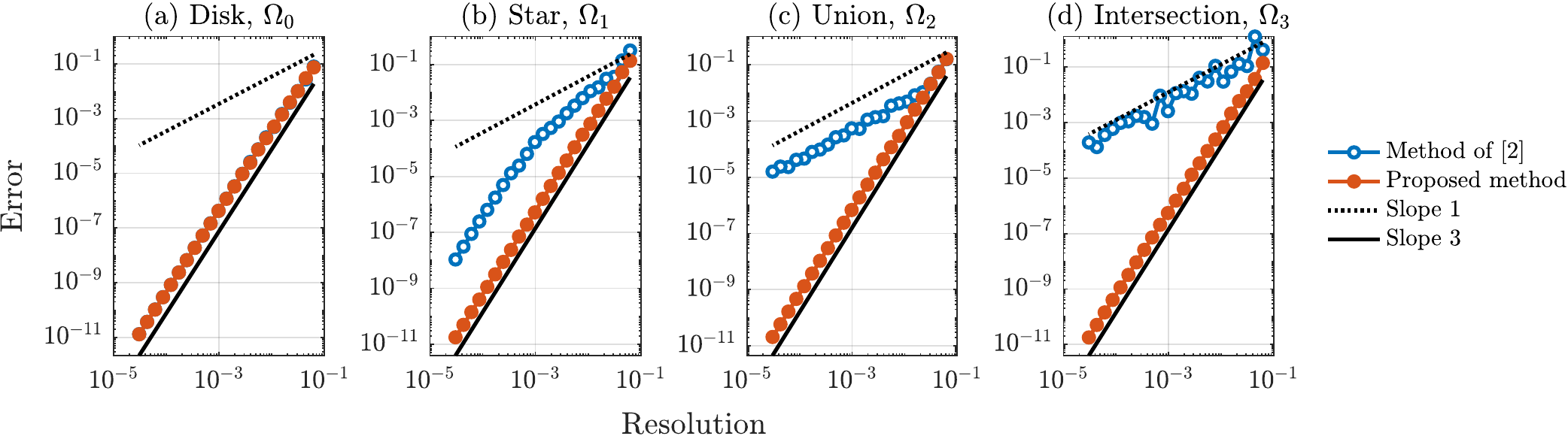}
\end{center}
\vspace{-.3cm}\vspace{-.3cm}\caption{\it Accuracy of the quadratic extrapolation \reviewerOne{(in the $L^\infty$ norm)} in two spatial dimensions measured in a narrow band of thickness $2\sqrt{\Delta x^2 + \Delta y^2}$ around an interface  using \reviewerOne{the approach of \cite{Aslam:04:A-partial-differenti} and the proposed approach}.} \label{fig::results::2d_quadratic}
\end{figure}

\begin{figure}[h]
\begin{center}
\subfigure[Disk, $\Omega_0$]{{\includegraphics[trim={0 0 2cm 0}, clip, width=.24\textwidth]{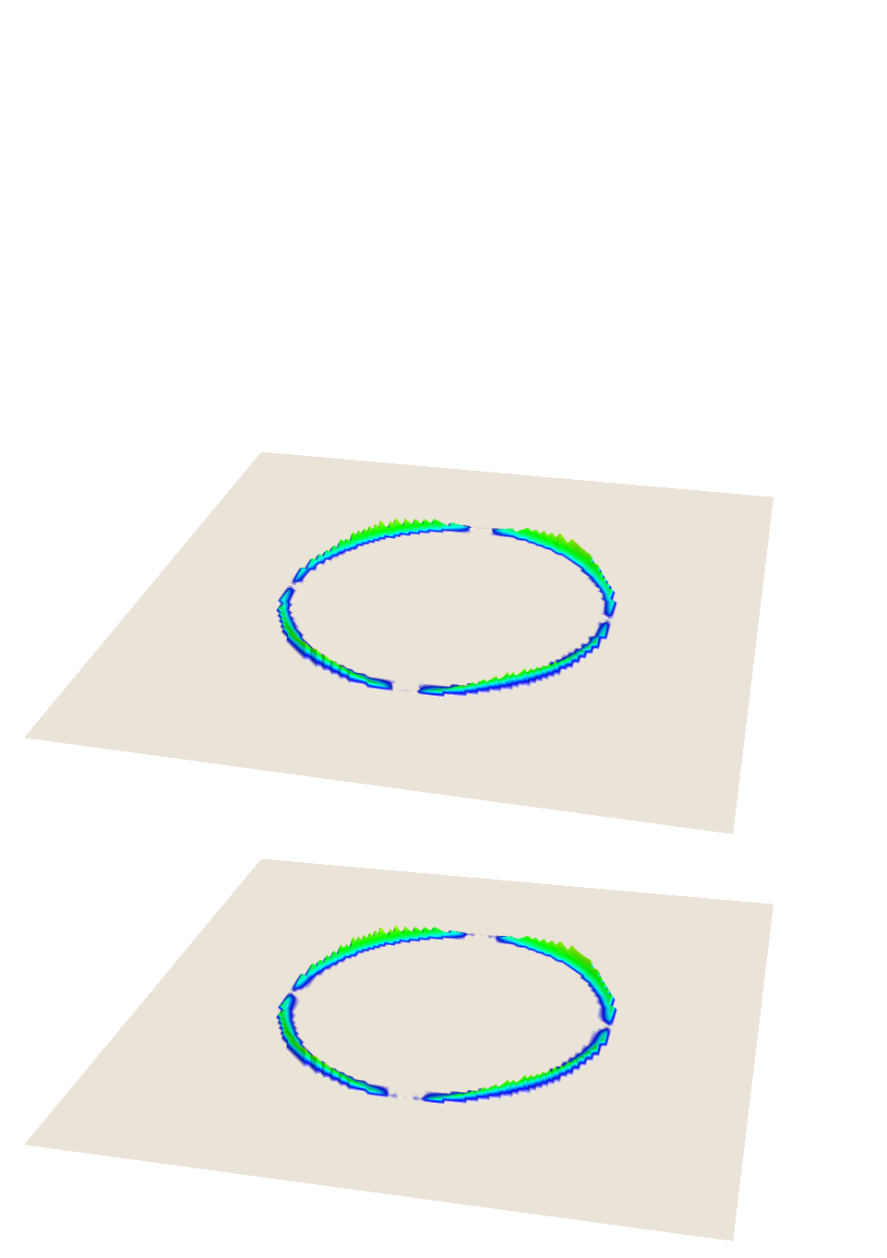}\label{fig:domain1_Bochkov}}} 
\subfigure[Star, $\Omega_1$]{{\includegraphics[trim={0 0 2cm 0}, clip, width=.24\textwidth]{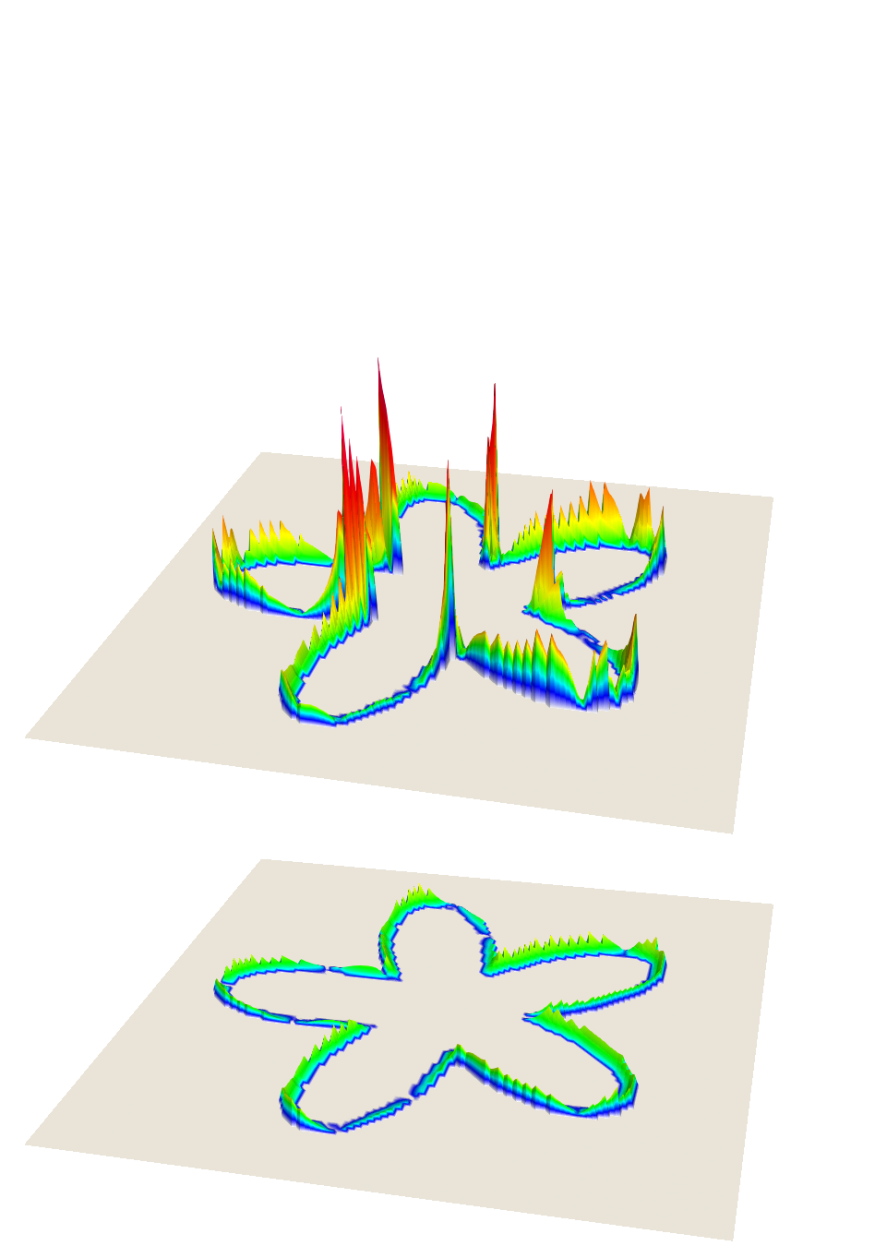}\label{fig:domain1_Bochkov}}} 
\subfigure[Union, $\Omega_2$]{{\includegraphics[trim={0 0 2cm 0}, clip, width=.24\textwidth]{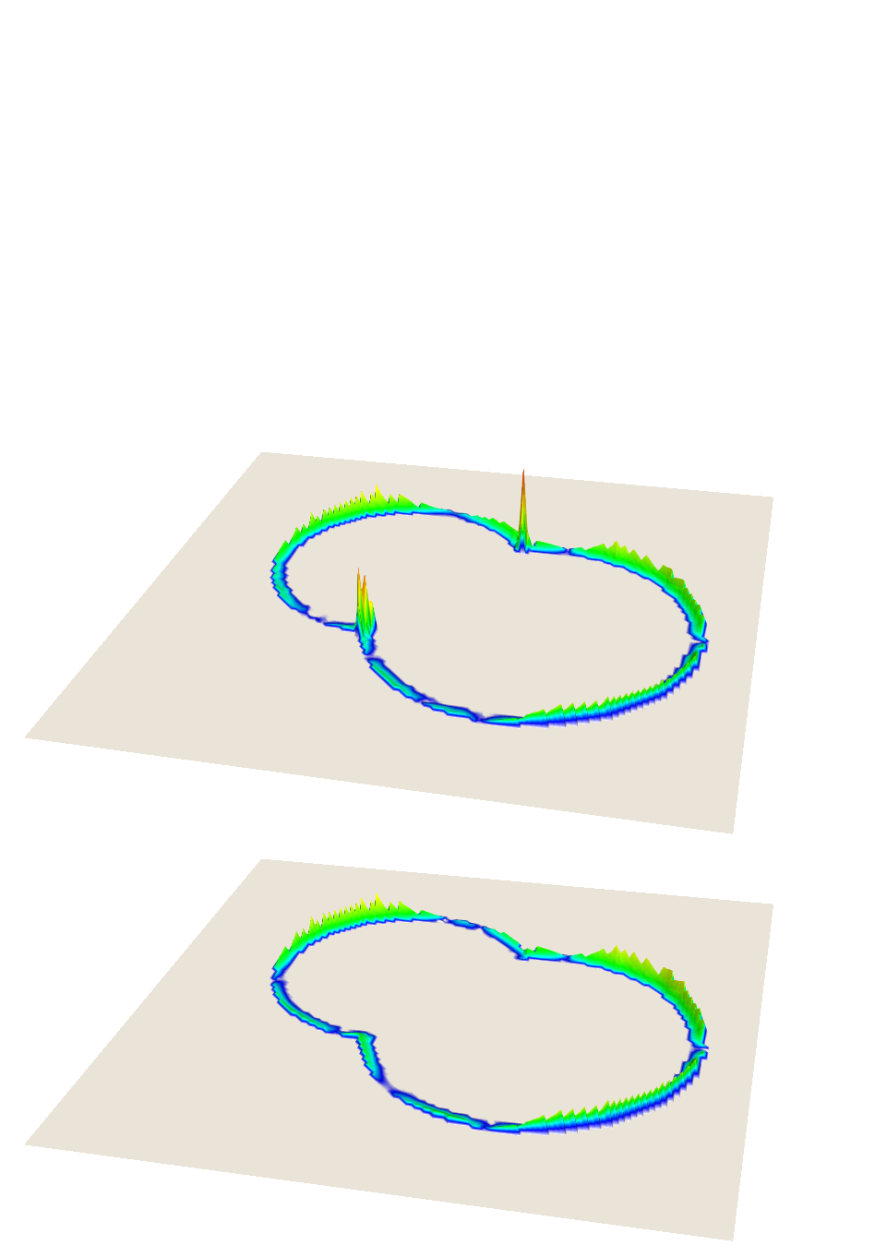}\label{fig:domain1_Bochkov}}} 
\subfigure[Intersection, $\Omega_3$]{{\includegraphics[trim={0 0 2cm 0}, clip, width=.24\textwidth]{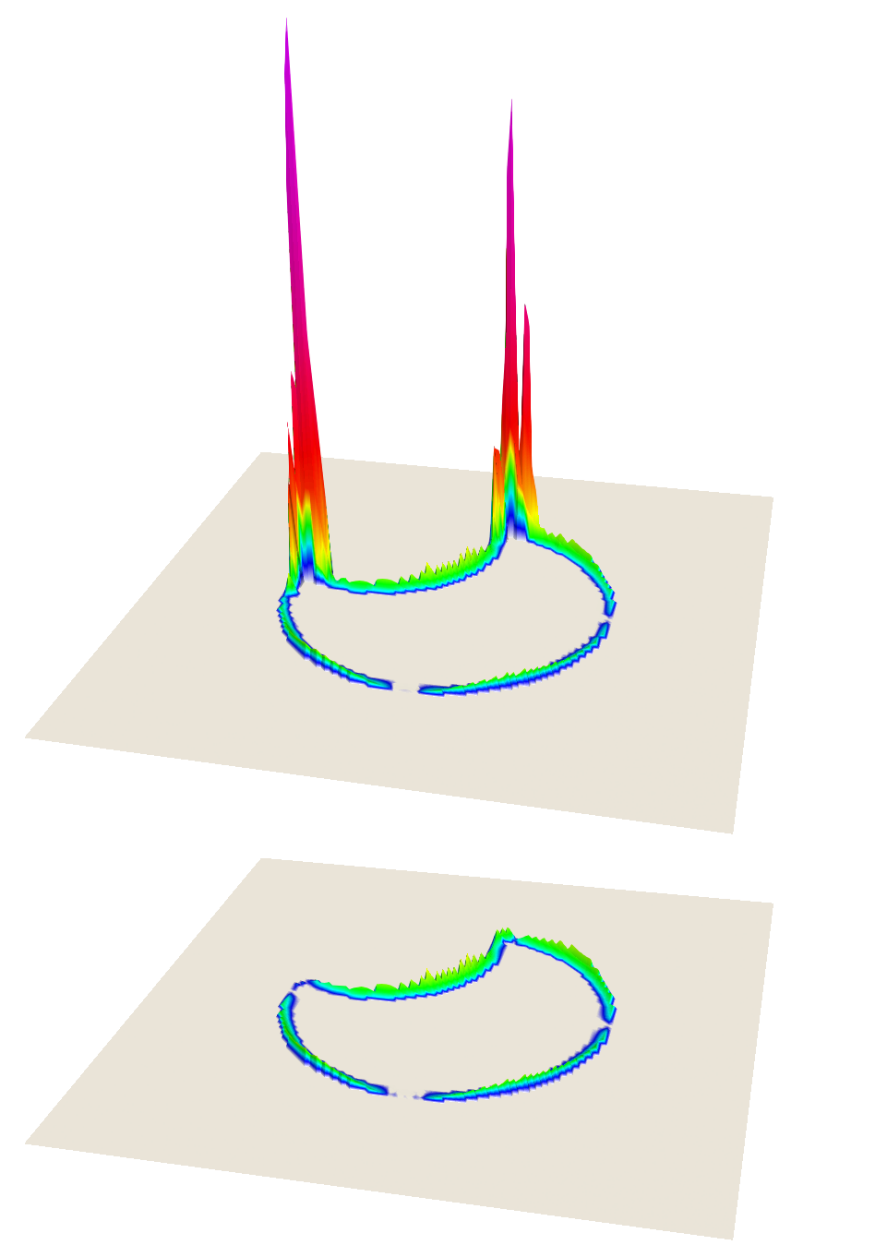}\label{fig:domain1_Bochkov}}}  
\end{center}
\vspace{-.3cm}\vspace{-.3cm}\caption{\it Comparison of error distributions in the case of the quadratic extrapolation on a $128^2$ grid. Top row: the approach of \cite{Aslam:04:A-partial-differenti}. Bottom row: the present approach. In each case the error is multiplied by a factor of 30 for visualization purpose.} \label{fig::Quadratic}
\end{figure}

In each case we consider a computational domain $\Omega=(-1, 1)\times(-1, 1)$. We define the function $q = \sin(\pi x)\cos(\pi y)$ inside every domain and extrapolate it in the outside region. Then the maximum difference between the exact values of $q$  and extrapolated ones\mylinelabel{rev1:Linf1}\reviewerOne{, that is, the $L^\infty$ norm of the error,} is computed within a band of thickness $2\sqrt{\Delta x^2 + \Delta y^2}$ in the outside regions. Figures \ref{fig::results::2d_linear} and \ref{fig::results::2d_quadratic} summarize the convergence behavior of the \mylinelabel{rev1:name9}\reviewerOne{ND-PDE} approach of \cite{Aslam:04:A-partial-differenti} and the proposed \mylinelabel{rev1:name10}\reviewerOne{WCD-PDE} approach. Figure \ref{fig::Quadratic} demonstrates the error distribution for both methods in the case of the quadratic extrapolation on a $128^2$ grid.

In case of the smooth domain $\Omega_0$ (disk), both approaches produce almost indistinguishable results attaining second- and third-order rates of convergence for the linear and quadratic extrapolation, respectively.

For the high-curvature domain $\Omega_1$ (star), both methods still reach optimal orders of convergence; however the \mylinelabel{rev1:name11}\reviewerOne{ND-PDE} approach demonstrates the optimal order of convergence only at relatively high grid resolutions when all geometric features are well-resolved. Moreover, for a given grid resolution the \mylinelabel{rev1:name12}\reviewerOne{WCD-PDE} approach produces results that are more than one order of magnitude more accurate in the case of the linear extrapolation and almost three orders of magnitude more accurate in the case of the quadratic extrapolation compared to the \mylinelabel{rev1:name13}\reviewerOne{ND-PDE} approach. Figure \ref{fig::Quadratic}b shows that the \mylinelabel{rev1:name14}\reviewerOne{ND-PDE} approach produces very large errors near regions with the highest curvature, while the error in the case of the \mylinelabel{rev1:name15}\reviewerOne{WCD-PDE} approach is much smaller and exhibits very little variation throughout all regions around the interface. 

The results are even more significantly improved with the proposed approach in the case of interfaces with kinks $\Omega_2$ (union) and $\Omega_3$ (intersection). Figures \ref{fig::Quadratic}c and \ref{fig::Quadratic}d show that the \mylinelabel{rev1:name16}\reviewerOne{ND-PDE} method of \cite{Aslam:04:A-partial-differenti} produces large errors near kinks; those errors are significantly reduced with the \mylinelabel{rev1:name17}\reviewerOne{WCD-PDE} approach. In particular, Figures \ref{fig::results::2d_linear}c-d and \ref{fig::results::2d_quadratic}c-d demonstrate that the second-order (third-order) accuracy of the linear (quadratic) extrapolations are recovered with the proposed approach; the rates of convergence for the approach of \cite{Aslam:04:A-partial-differenti} are close to first order, which corresponds to the constant extrapolation, due to the fact that errors near kinks do not decrease despite grid refinement \mylinelabel{rev2:first-order}\reviewerTwo{and the apparent first order of convergence is only because  the neighborhood in which errors are computed is shrinking closer to the domain}.

\section{Numerical Results in Three Spatial Dimensions} \label{sec:Numerical_Examples_3d}

\begin{figure}[h]
\begin{center}
\subfigure[Sphere, $\tilde{\Omega}_0$]{{\includegraphics[width=.24\textwidth]{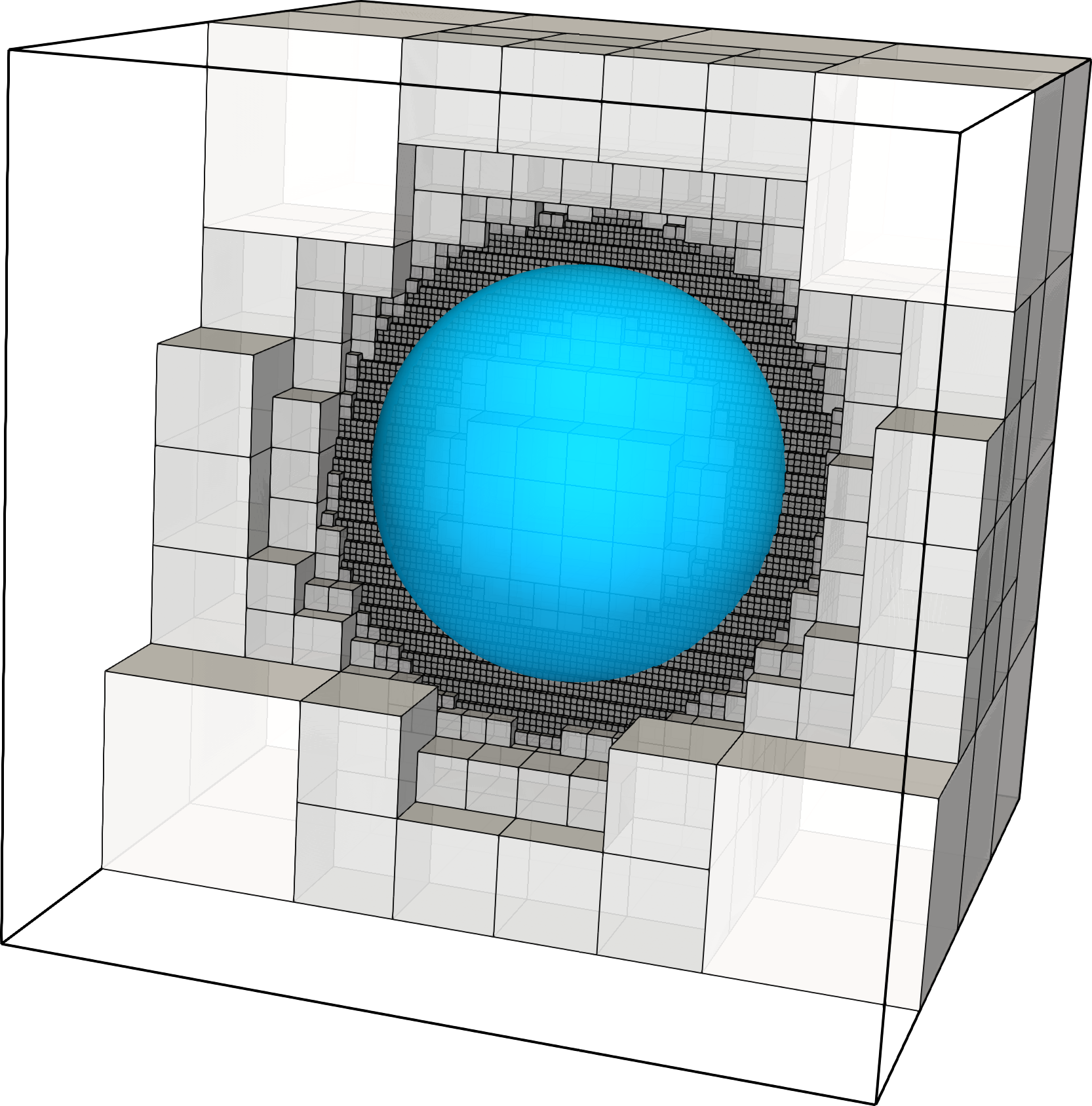}\label{fig:octree_sphere}}} 
\subfigure[Star, $\tilde{\Omega}_1$]{{\includegraphics[width=.24\textwidth]{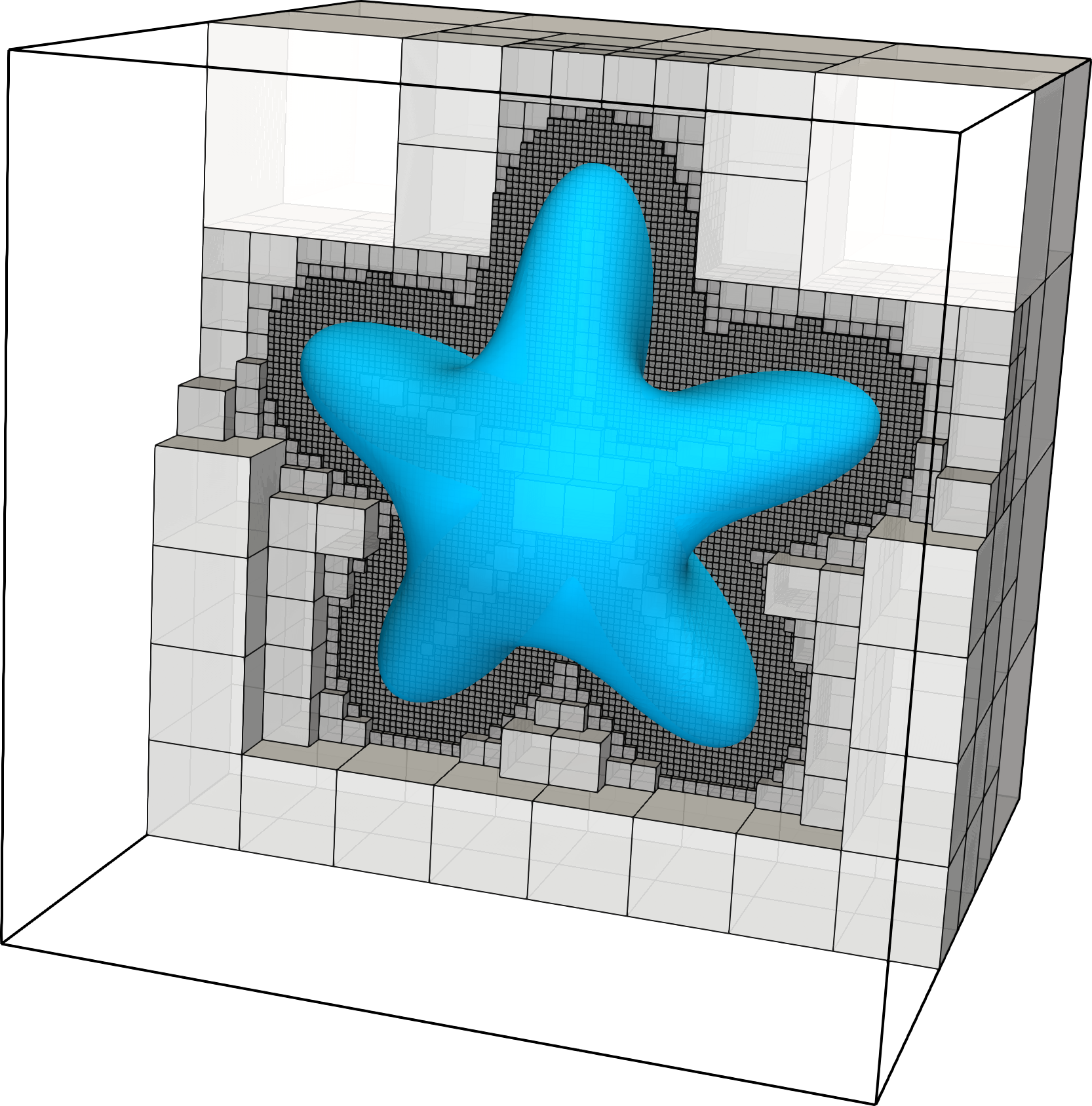}\label{fig:octree_star}}} 
\subfigure[Union, $\tilde{\Omega}_2$]{{\includegraphics[width=.24\textwidth]{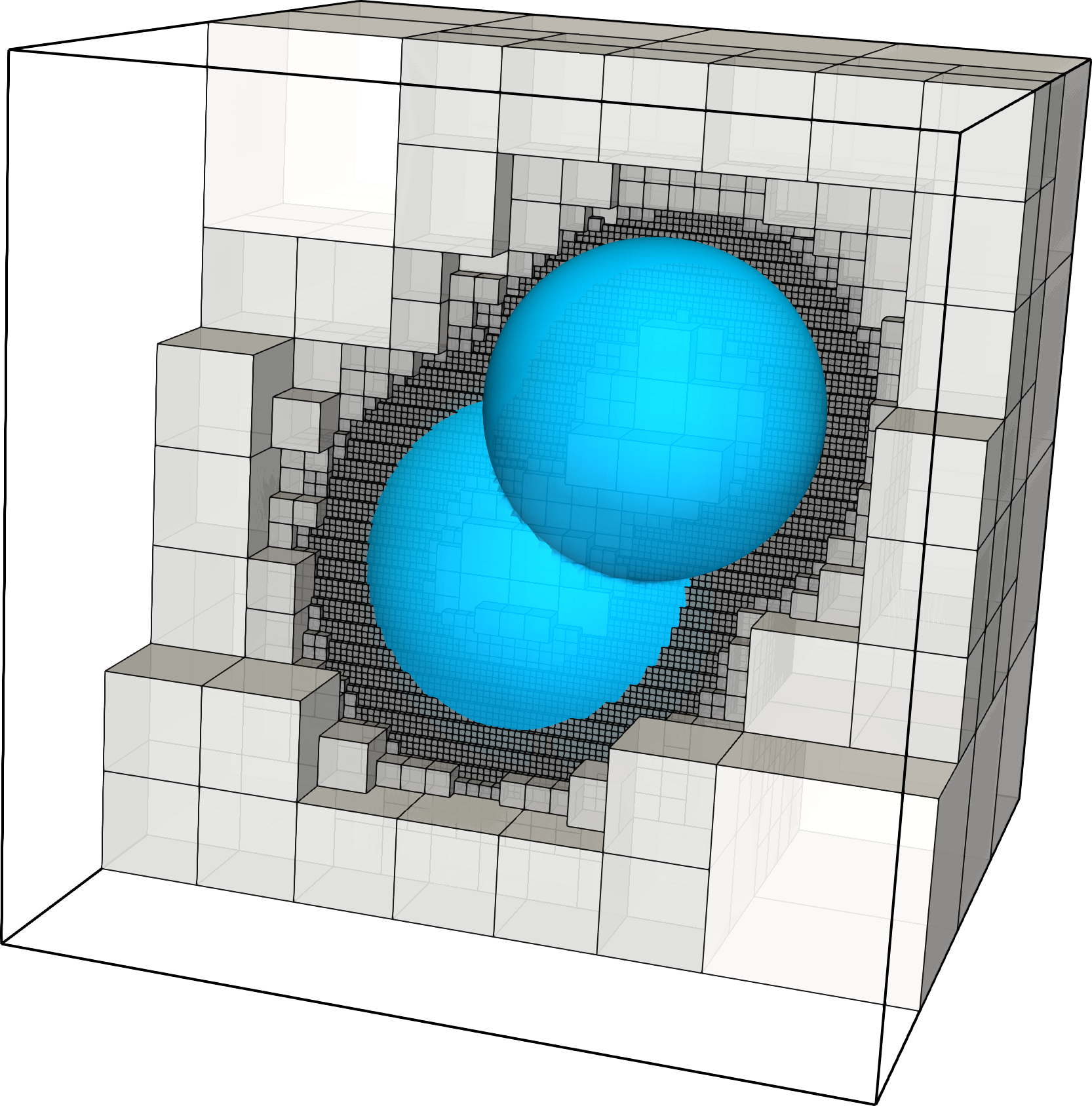}\label{fig:octree_union}}} 
\subfigure[Intersection, $\tilde{\Omega}_3$]{{\includegraphics[width=.24\textwidth]{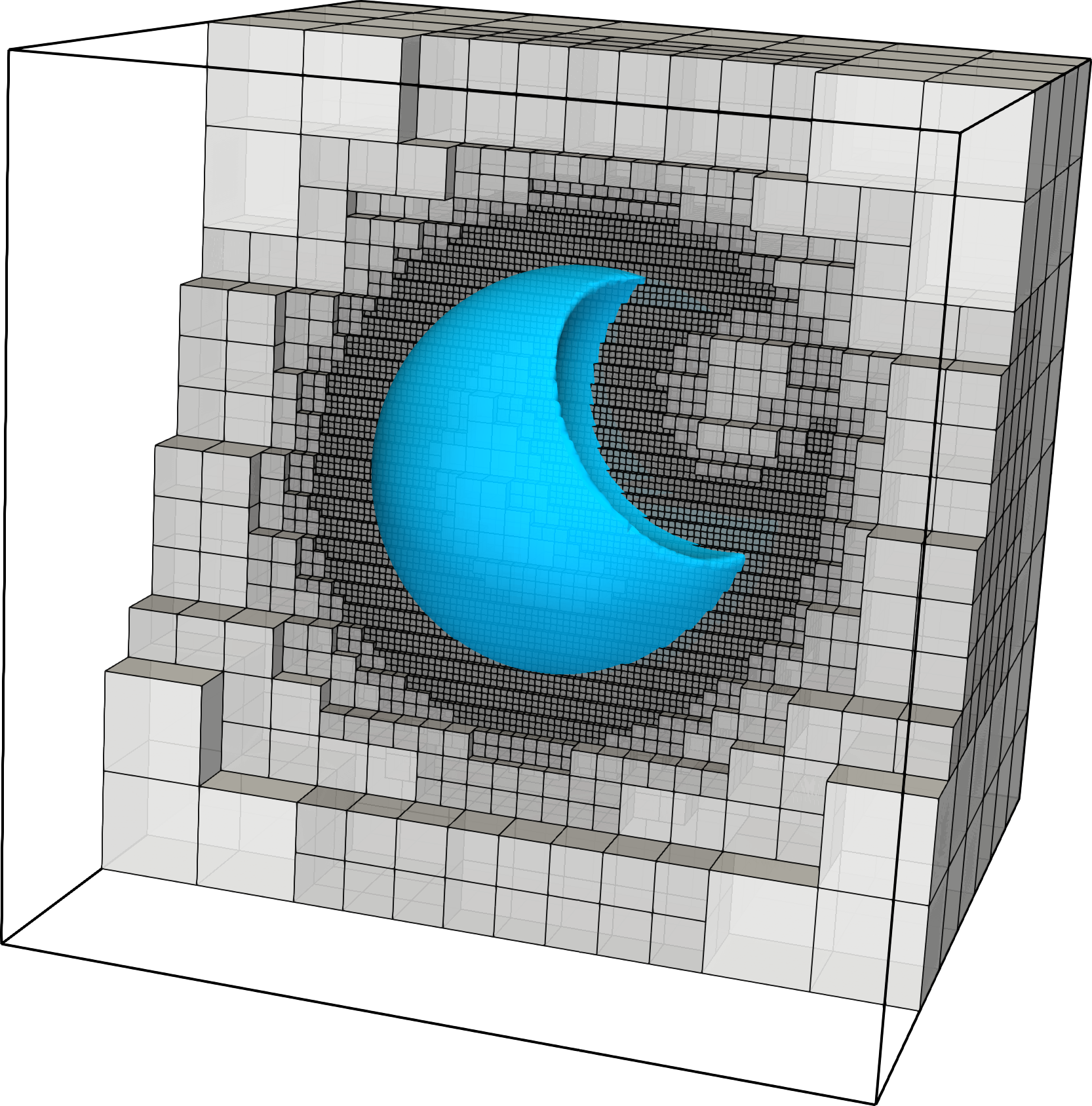}\label{fig:octree_intersection}}} 
\end{center}
\vspace{-.3cm}\vspace{-.3cm}\caption{\it Irregular domains considered in section \ref{sec:Numerical_Examples_3d}
 along with the octree grids refined near their boundaries.} \label{fig::octree_domains}
\end{figure}

We consider three different domains, $\tilde{\Omega}_1$, $\tilde{\Omega}_2$ and $\tilde{\Omega}_3$, that present high-curvature features or kinks in three spatial dimensions. In addition, we consider a smooth spherical domain $\tilde{\Omega}_0$ with center $(0,0,0)$ and radius $0.501$. The definition of those domains are given by the level-set functions:
\begin{linenomath*}
\begin{align*}
\begin{aligned}
\tilde{\phi}_0(x, y, z) &= \sqrt{x^2 + y^2 + z^2}-0.501, \\
\tilde{\phi}_1(x, y, z) &= \sqrt{x^2 + y^2 + z^2}-0.501 - 0.15 \frac{y^5 + 5x^4y-10x^2y^3}{\left(x^2 + y^2 + z^2\right)^{\frac{5}{2}}} \cos \left(\frac{\pi}{2}\frac{z}{0.501} \right), \\
\tilde{\phi}_2(x, y, z) &= \min\left(\sqrt{(x+.1)^2 + (y+.3)^2 + (z+.2)^2}-0.501,
                   \sqrt{(x-.2)^2 + (y-.2)^2 + (z-.1)^2}-0.401\right), \\
\tilde{\phi}_3(x, y, z) &= \max\left(\sqrt{x^2 + y^2 + z^2}-0.501,
                   \sqrt{(x-.4)^2 + (y-.3)^2 + (z-.2)^2}-0.401\right),
\end{aligned}
\end{align*}
\end{linenomath*}
Figure \ref{fig::octree_domains} depicts those domains along with the octree grid refined near their boundaries.

Similar to the two-dimensional examples, we consider a computational domain $\Omega=(-1, 1)^3$. We extrapolate the function $q = \sin(\pi x)\cos(\pi y) \exp(z)$ from the inside to the outside for every domain and compute the difference between the exact values of $q$ and the extrapolated ones\mylinelabel{rev1:Linf2}\reviewerOne{, that is, the $L^\infty$ norm of the error,} within a band of thickness $2\sqrt{\Delta x^2 + \Delta y^2 + \Delta z^2}$ in the outside region.

\begin{figure}[!h]
\begin{center}
\includegraphics[width=.99\textwidth]{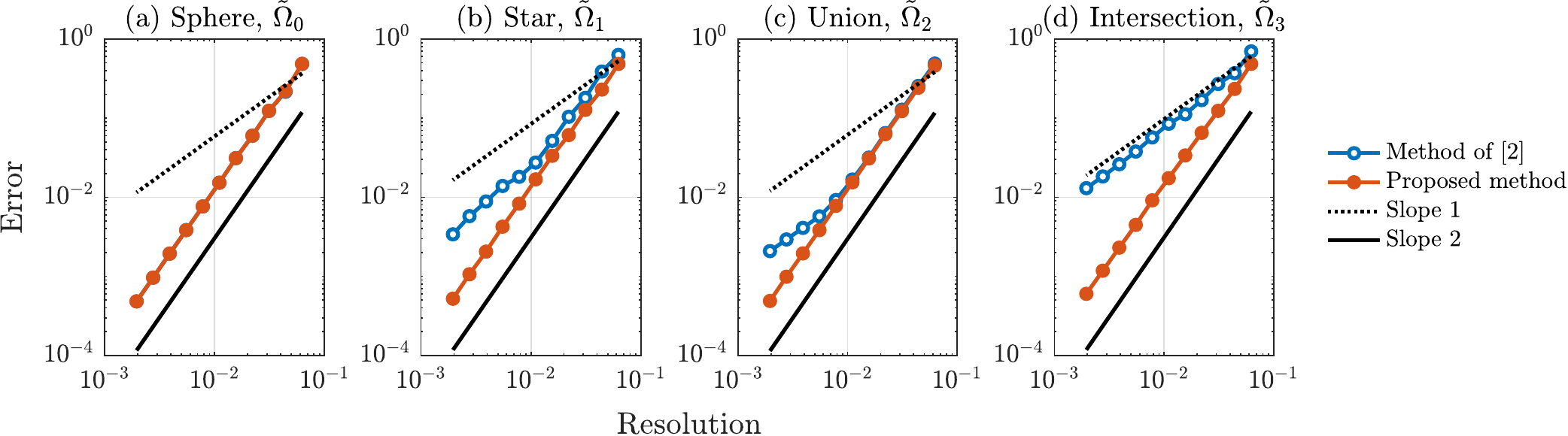}
\end{center}
\vspace{-.3cm}\vspace{-.3cm}\caption{\it Accuracy of the linear extrapolation \reviewerOne{(in the $L^\infty$ norm)} in three spatial dimensions measured in a narrow band of thickness $2\sqrt{\Delta x^2 + \Delta y^2 + \Delta z^2}$ around an interface  using \reviewerOne{the approach of \cite{Aslam:04:A-partial-differenti} and the proposed approach}.} \label{fig::results::3d_linear}
\end{figure}

\begin{figure}[!h]
\begin{center}
\includegraphics[width=.99\textwidth]{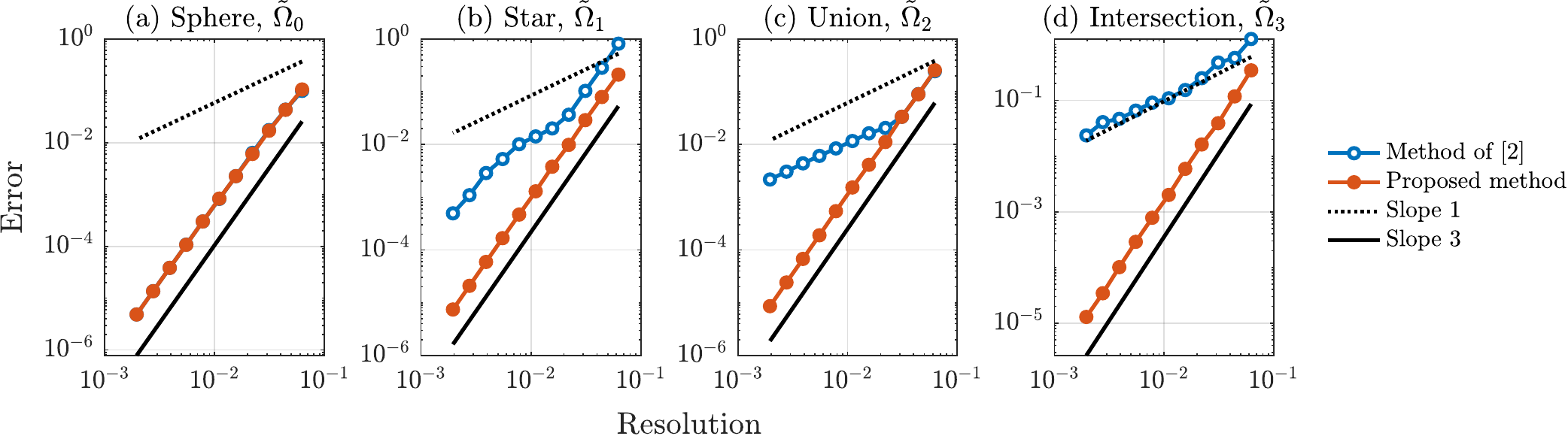}
\end{center}
\vspace{-.3cm}\vspace{-.3cm}\caption{\it Accuracy of the quadratic extrapolation \reviewerOne{(in the $L^\infty$ norm)} in three spatial dimensions measured in a narrow band of thickness $2\sqrt{\Delta x^2 + \Delta y^2 + \Delta z^2}$ around an interface  using \reviewerOne{the approach of \cite{Aslam:04:A-partial-differenti} and the proposed approach}.} \label{fig::results::3d_quadratic}
\end{figure}

Conclusions similar to the two dimensional case can be drawn from the results in figures \ref{fig::results::3d_linear} and \ref{fig::results::3d_quadratic}. Specifically, for a smooth and well-resolved domain $\tilde{\Omega}_0$ (sphere) both approaches produce almost indistinguishable results with optimal order of convergence (second and third for the linear and quadratic extrapolations, respectively). When the interface curvature is high ($\tilde{\Omega}_1$, star) the \mylinelabel{rev1:name18}\reviewerOne{WCD-PDE} approach produce extrapolated fields that are several orders of magnitude more accurate than for the \mylinelabel{rev1:name19}\reviewerOne{ND-PDE} approach. For geometries with sharp features $\tilde{\Omega}_2$ (union) and $\tilde{\Omega}_3$ (intersection) only the \mylinelabel{rev1:name20}\reviewerOne{WCD-PDE} approach demonstrates optimal orders of convergence, while for the \mylinelabel{rev1:name21}\reviewerOne{ND-PDE} approach the rate of convergence is stuck to 1.

\reviewerOne{
\section{Application to Solving the Diffusion Equation in Time-Dependent Domains} 
\label{sec:Numerical_Diffusion}}

\begin{figure}[h]
\begin{center}
\subfigure[Smooth moving object]{{\includegraphics[width=.3\textwidth]{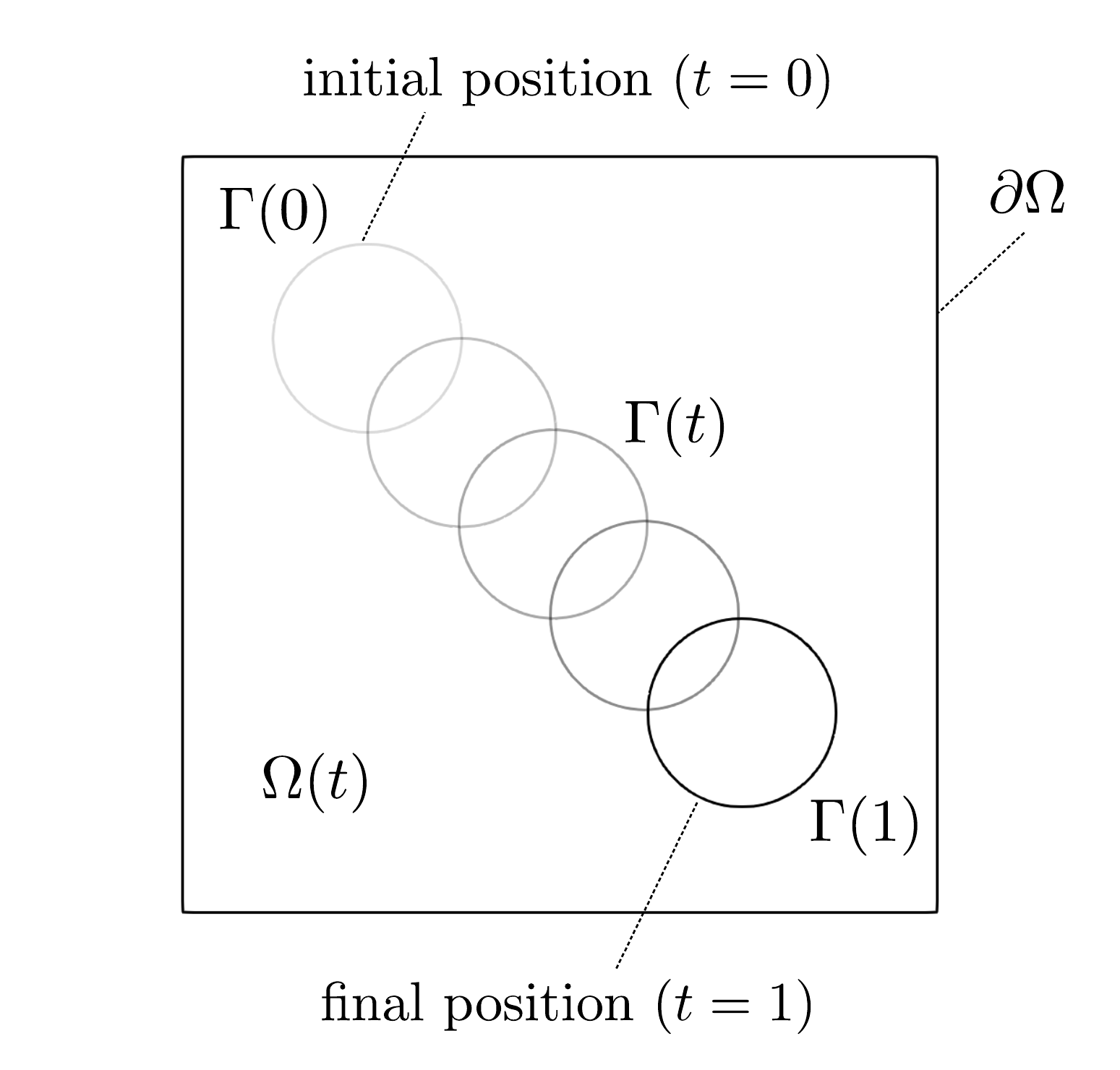}\label{fig:octree_sphere}}} 
\subfigure[Moving object with corners]{{\includegraphics[width=.3\textwidth]{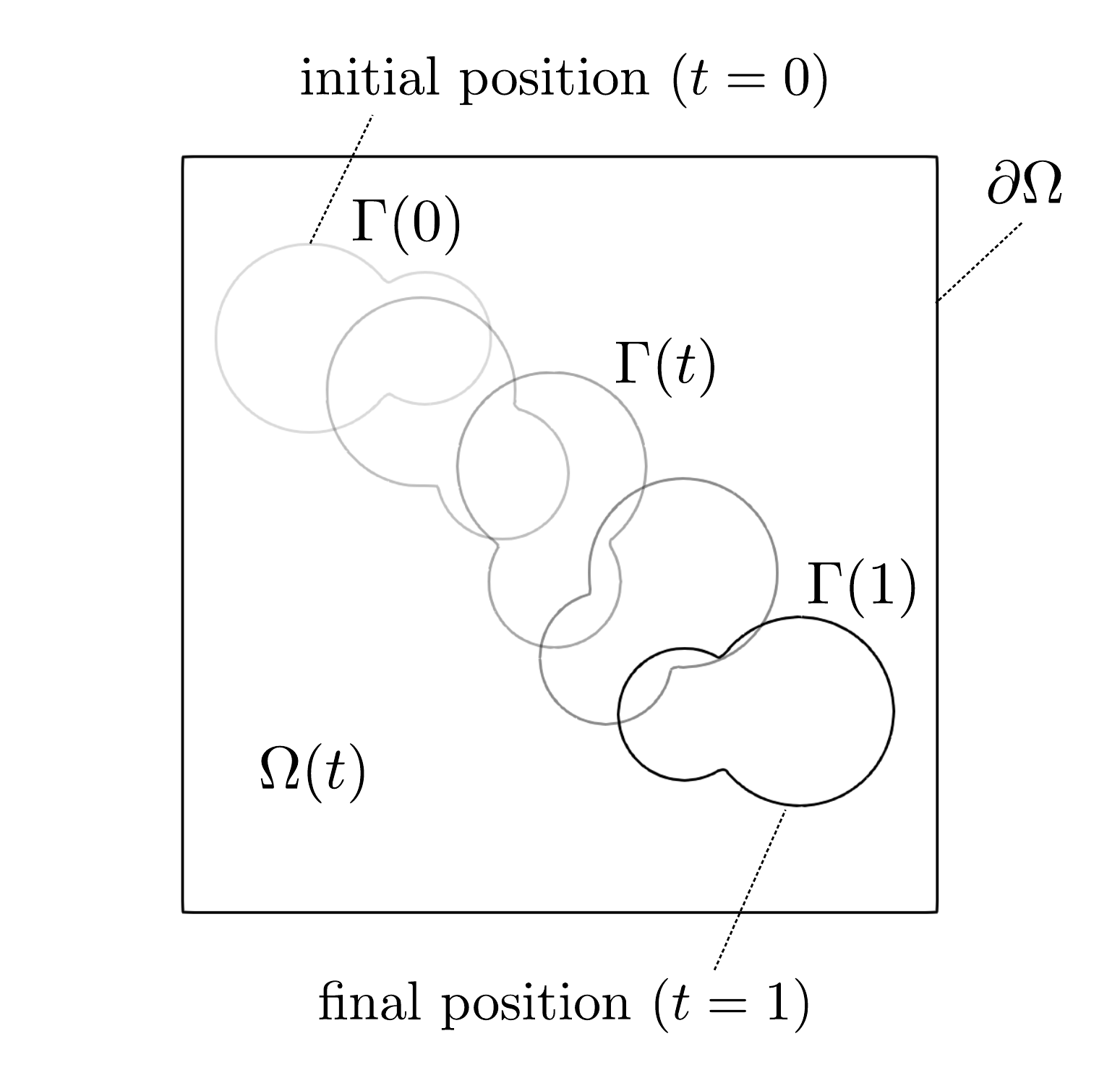}\label{fig:octree_star}}} 
\end{center}
\vspace{-.3cm}\vspace{-.3cm}\caption{\it Problem geometry in diffusion equation example from Section \ref{sec:Numerical_Diffusion}.} \label{fig::diffusion_geometry}
\end{figure}

\reviewerOne{
In order to illustrate the importance of accurate extrapolation near sharp corners in moving interface problems, we present a simple example of solving diffusion equation around a moving object that may have a non-smooth boundary. Specifically, we consider a two-dimensional rectangular region $[-1;1] \times [-1;1]$ 
and an object that moves diagonally from its starting position at $(x_s, y_s) = (-0.51, 0.52)$ at time $t=0$ to the final position $(x_f, y_f) = (0.49, -0.48)$ at time $t=1$ while making half a turn around its center as demonstrated in Figure \ref{fig::diffusion_geometry}. A diffusion equation subject to Neumann boundary conditions is solved in the rectangular box excluding the region occupied by the moving object. While the problem at hand does not correspond to any specific practical application, it represents a prototypical situation arising in simulation of more relevant (and more complex) processes and at the same time allows a precise analysis of numerical errors. Note that a non-deformable shape is considered for the sake of simplicity, we expect the main conclusions to hold true in more general cases, e.g. multiphase flow with triple junction points.
}

\reviewerOne{
The Eulerian framework is employed, more precisely, the region $[-1;1] \times [-1;1]$ is discretized into a static uniform rectangular grid with $N$ nodes in each Cartesian direction while the object is implicitly described by a time-dependent level-set function.
}

\reviewerOne{
Suppose the shape of the moving object in its local system of coordinate $\vect{\xi}$ is described by a level-set function $\phi_0(\vect{\xi})$ (such that $\phi_0(\vect{\xi}) > 0$ inside the object). Then the object's motion in the global system of coordinates $\vect{r}$ can be expressed by a time-dependent level-set function $\phi(t, \vect{r}) = \phi_0 (\vect{\xi}(t, \vect{r}))$, where global-to-local coordinate transformation $\vect{\xi}(t, \vect{r})$ is given by:
\begin{linenomath*}
\begin{align*}  
  \vect{\xi} (t,\vect{r}) = 
  \begin{pmatrix}
    \xi (t,\vect{r}) \\ \eta (t,\vect{r})
  \end{pmatrix}
  = 
  \begin{pmatrix}
    x - (x_s + t(x_f-x_s)) \\
    y - (y_s + t(y_f-y_s)) 
  \end{pmatrix}
  \begin{pmatrix}
    \cos (\pi t) & \sin (\pi t) \\
    -\sin (\pi t) & \cos (\pi t)
  \end{pmatrix}.
\end{align*}
\end{linenomath*}
The solution domain $\Omega(t)$ can then be defined as:
\begin{linenomath*}
\begin{align*}
  \Omega(t) = \left\{ \vect{r} \in [-1;1]\times[-1;1] \, : \, \phi(t, \vect{r}) < 0 \right\}.
\end{align*}
\end{linenomath*}
The boundary of the computational box is denoted as $\partial\Omega$ and the boundary of the moving object is denoted as $\Gamma(t)$.
}

\reviewerOne{
In order to investigate the influence of non-smooth interface, we consider two choices of moving object, one having a smooth boundary, a disk of radius $0.25$, and another one having a non-smooth boundary, a union of two disk with radii $0.25$ and $0.175$, motivated by multimaterial compound bubbles. In the first case the level-set function of the object is given by (in the local system of coordinates):
\begin{linenomath*}
\begin{align*}
  \phi_0^{\textrm{smooth}} (\vect{\xi})  = (\vect{\xi}) = r_0 - \sqrt{\xi^2 + \eta^2},
\end{align*}
\end{linenomath*}
while in the latter case:
\begin{linenomath*}
\begin{align*}
  \phi_0^{\textrm{non-smooth}} (\vect{\xi}) = \max \left(
  r_0 - \sqrt{\left(\xi + \xi_0 \right)^2 + \eta^2}, 
  r_0q - \sqrt{\left(\xi - \xi_0 \right)^2 + \eta^2} \right),
\end{align*}
\end{linenomath*}
where $r_0 = 0.25$, $ \xi_0 = \hf r_0 \sqrt{1+q^2}$ and $q = 0.7$.
}

\reviewerOne{
We choose the following test solution:
\begin{linenomath*}
\begin{align*}
  u(t,\vect{r}) = \sum_{i=0}^{3} \sum_{j=0}^{3} a_{i,j} \cos \left( \frac{i \pi}{2} (x+1) \right) \cos \left( \frac{j \pi}{2} (y+1) \right) \exp \left( -D\frac{\pi^2}{4} (i^2+j^2) t \right),
\end{align*}
\end{linenomath*}
with 
\begin{linenomath*}
\begin{align*}
\begin{pmatrix}
a_{0,0} & a_{0,1} & a_{0,2} & a_{0,3} \\
a_{1,0} & a_{1,1} & a_{1,2} & a_{1,3} \\
a_{2,0} & a_{2,1} & a_{2,2} & a_{2,3} \\
a_{3,0} & a_{3,1} & a_{3,2} & a_{3,3}
\end{pmatrix}
=
\begin{pmatrix}
-0.5 & -0.1 & -0.5 &  0.6 \\ 
-0.6 & -0.5 & -0.1 & -0.1 \\ 
 0.2 & -0.2 & -0.2 &  0.4 \\ 
 0.1 & -0.9 &  0.8 &  0.4
\end{pmatrix},
\end{align*}
\end{linenomath*}
which satisfies a homogeneous diffusion equation:
\begin{linenomath*}
\begin{alignat}{2}
\label{eq:diff_eq}
  \ddt{u} &= D \lap u,&& \quad \textrm{ for } t \in [0,1], \, \vect{r} \in \Omega(t),
\end{alignat}
\end{linenomath*}
with initial conditions:
\begin{linenomath*}
\begin{align*}
  u(0, \vect{r}) = \sum_{i=0}^{3} \sum_{j=0}^{3} a_{i,j} \cos \left( \frac{i \pi}{2} (x+1) \right) \cos \left( \frac{j \pi}{2} (y+1) \right), \quad \textrm{ for } \vect{r} \in  \partial \Omega(0)
\end{align*}
\end{linenomath*}
and boundary conditions:
\begin{linenomath*}
\begin{alignat*}{2}
  D \ddn{u} &= 0,&& \quad \textrm{ for } t \in [0,1], \, \vect{r} \in  \partial \Omega,\\
  D \ddn{u} &= g(t,\vect{r}),&& \quad \textrm{ for } t \in [0,1], \, \vect{r} \in  \Gamma(t),
\end{alignat*}
\end{linenomath*}
where function $g(t,\vect{r})$ is given by:
\begin{linenomath*}
\begin{align*}
  g(t, \vect{r}) = D\frac{\nabla \phi(t,\vect{r}) \cdot \nabla u(t,\vect{r})}{\left| \nabla \phi(t,\vect{r}) \right|}.
\end{align*}
\end{linenomath*}
}

\reviewerOne{
The time range $[0;1]$ is discretized into time layers $t_n$, $n = 0,1,2\ldots$, where the time step between adjacent time layers $\Delta t_{n+1} = t_{n+1} - t_{n}$ is determined such that the maximum displacement of the moving object boundary during the given time step is expected not to exceed a user-defined fraction $f$ of the grid cell diagonal $\sqrt{\Delta x^2 + \Delta y^2}$, that is:
\begin{linenomath*}
\begin{align*}
  \Delta t_{n+1} = \frac{f \sqrt{\Delta x^2 + \Delta y^2}}{\max\limits_{\Gamma(t_{n})} v_{\vect{n}}(t_{n}, \vect{r})},
\end{align*}
\end{linenomath*}
where $v_{\vect{n}}$ denotes the normal velocity of the object's boundary. Specifically, in this example $f=0.8$ is taken.
}

\reviewerOne{
We use the second-order variable-step backward differentiation formula (BDF2) for discretizing the diffusion equation \eqref{eq:diff_eq} in time, that is:
\begin{linenomath*}
\begin{align}
\label{eq:diff_eq_discrete}
  \left( \alpha_0 - \Delta t_n \lap \right) u^{n} = -\alpha_1 u^{n-1} -\alpha_2 u^{n-2},
\end{align}
\end{linenomath*}
where 
\begin{linenomath*}
\begin{align*}
  \alpha_0 = \frac{1+2\rho}{1+\rho}, \quad \alpha_1 = -(1+\rho), \quad \alpha_2 = \frac{\rho^2}{1+\rho} \quad \textrm{and} \quad \rho = \frac{\Delta t_n}{\Delta t_{n-1}},
\end{align*}
\end{linenomath*}
and use the superconvergent second-order accurate method of \cite{Bochkov;Gibou:19:Solving-the-Poisson-} (which is specifically designed to handle irregular domains with non-smooth boundaries) for solving the resulting Poisson-type equation \eqref{eq:diff_eq_discrete}.
}

\reviewerOne{
Solution of equation \eqref{eq:diff_eq_discrete} produces values of $u^n$ at all grid nodes that belong to the current solution domain $\Omega(t_n)$. However, as the object moves in time some of grid nodes outside of $\Omega(t_n)$ may become part of $\Omega(t_{n+1})$ or $\Omega(t_{n+2})$ and solving at subsequent time layers $t_{n+1}$ and $t_{n+2}$ would require valid values of $u^n$ at such grid nodes. This is typically addressed in free boundary value problems by smoothly extrapolating $u^n$ into some neighborhood of $\Omega(t_n)$. In this example we quadratically extrapolate solutions using both the proposed in this work WCD-PDE approach and the ND-PDE approach of \cite{Aslam:04:A-partial-differenti}. Also, since solving Poisson-type equations on irregular domains produces additional errors, we generate a reference solution where instead of performing extrapolation of numerical values we fill the grid nodes outside of the current solution domain $\Omega(t_n)$ with exact values given by the analytical solution.
}

\reviewerOne{
We investigate the influence of extrapolation procedures on the accuracy of the numerical solution and its gradient. In order not to measure the error of extrapolation procedure itself but rather only its influence on solving the diffusion equation, the solution error at time layer $t_n$ is calculated only for grid nodes in $\Omega(t_n)$ and the gradient is calculated only using values from $\Omega(t_n)$ as well.
}

\begin{figure}[!h]
\begin{center}
\subfigure[Using exact values]{{\includegraphics[trim={10cm 17cm 10cm 25cm}, clip, width=.25\textwidth]{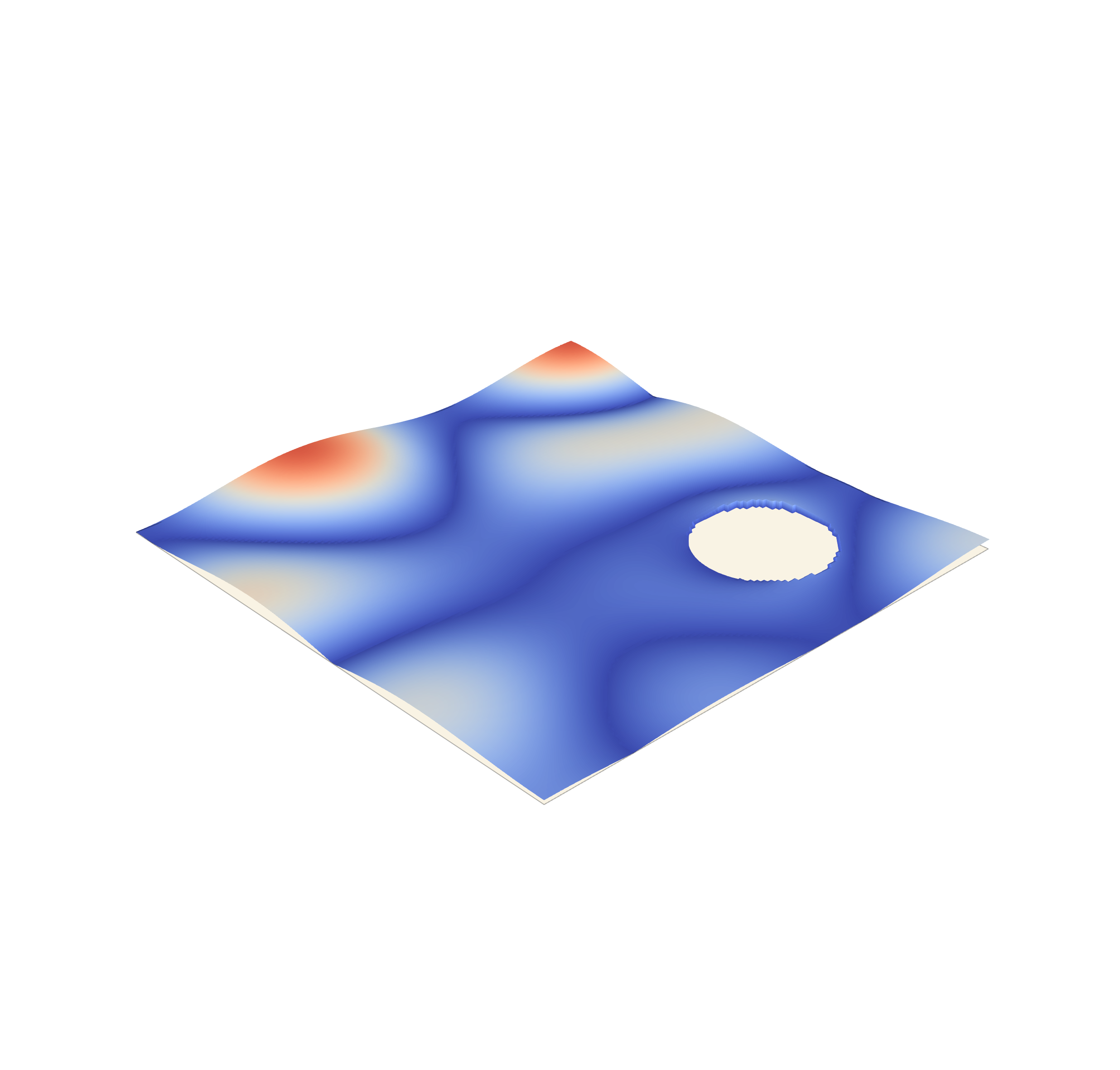}}} 
\subfigure[Using method of \cite{Aslam:04:A-partial-differenti}]{{\includegraphics[trim={10cm 17cm 10cm 25cm}, clip, width=.25\textwidth]{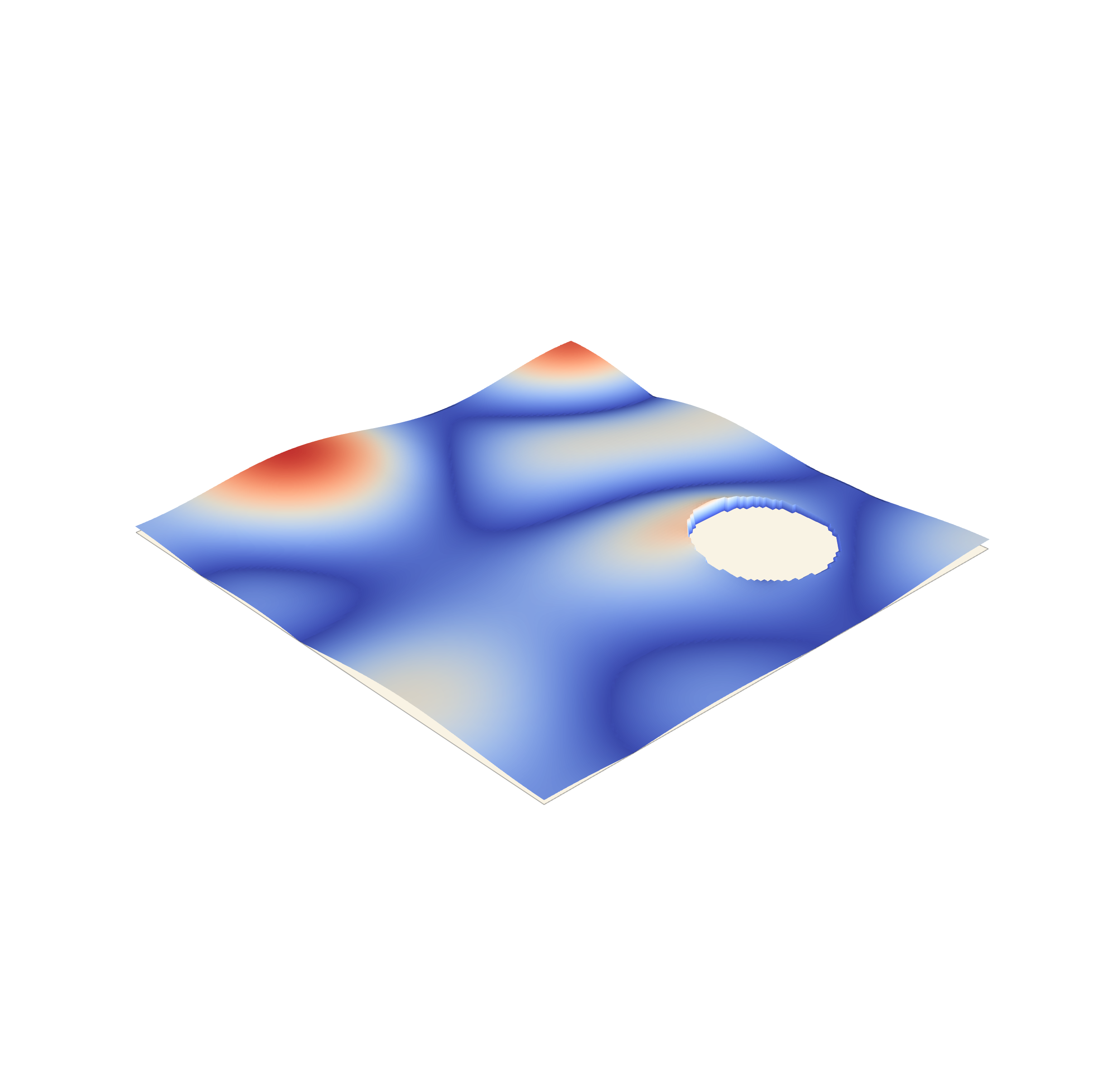}}} 
\subfigure[Using proposed method]{{\includegraphics[trim={10cm 17cm 10cm 25cm}, clip, width=.25\textwidth]{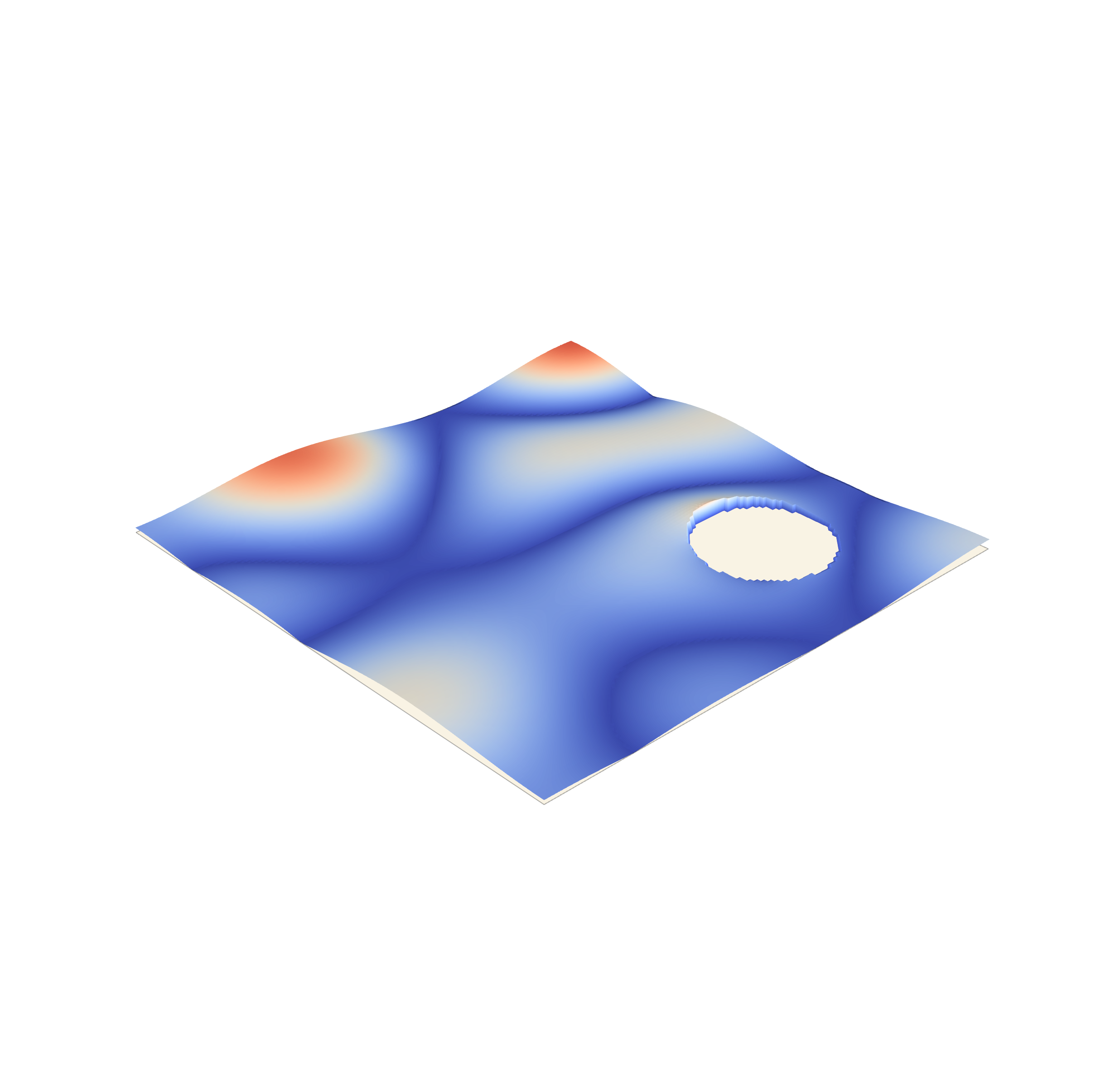}}} 
\end{center}
\vspace{-.3cm}\vspace{-.3cm}\caption{\it Error distribution at the final time moment ($t=1$) in case of a smooth moving object using different extrapolation approaches.} \label{fig::diffusion_circle_error}
\end{figure}

\begin{figure}[!h]
\begin{center}
\subfigure[Using exact values]{{\includegraphics[trim={10cm 17cm 10cm 6cm}, clip, width=.25\textwidth]{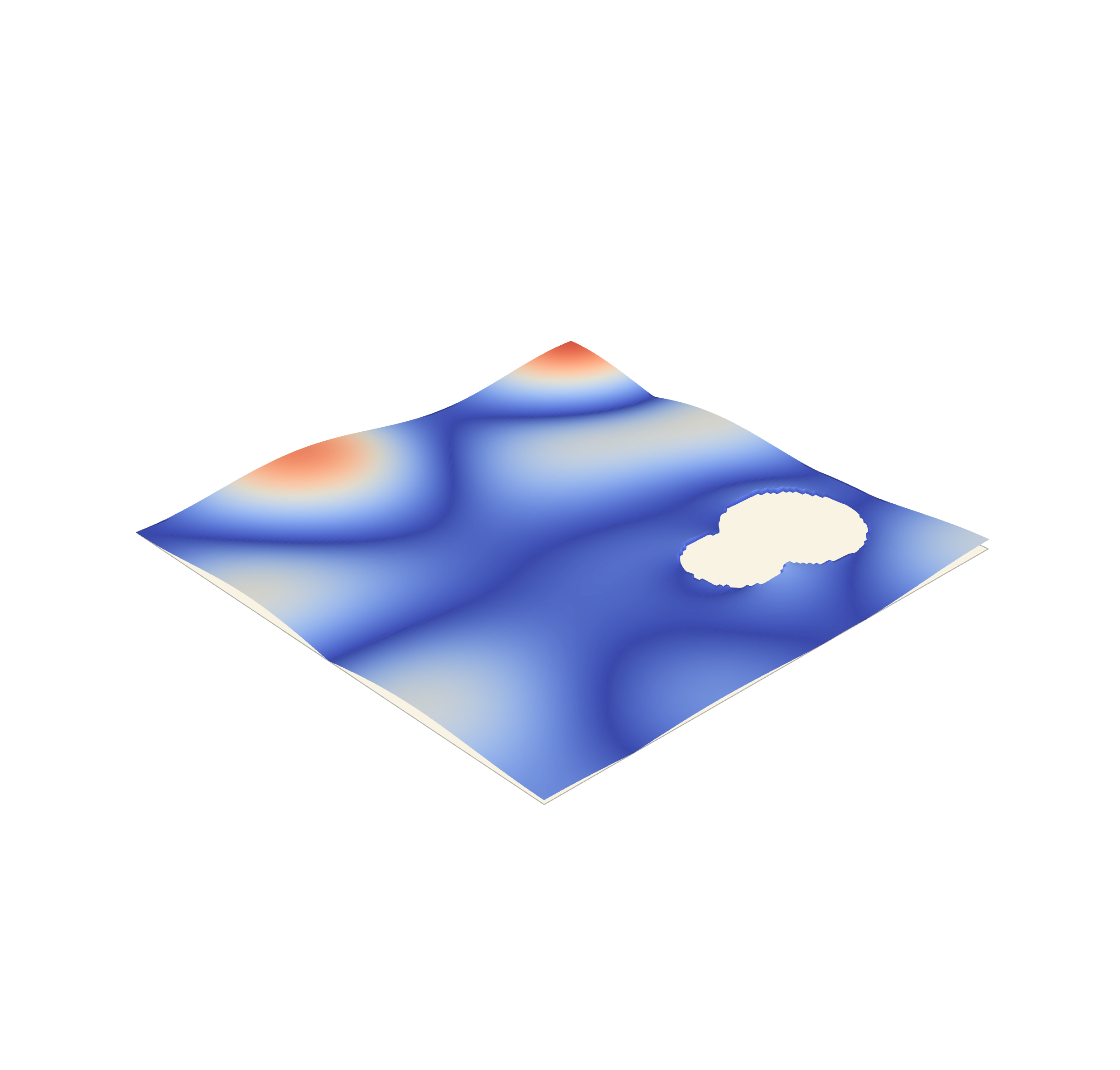}}} 
\subfigure[Using method of \cite{Aslam:04:A-partial-differenti}]{{\includegraphics[trim={10cm 17cm 10cm 6cm}, clip, width=.25\textwidth]{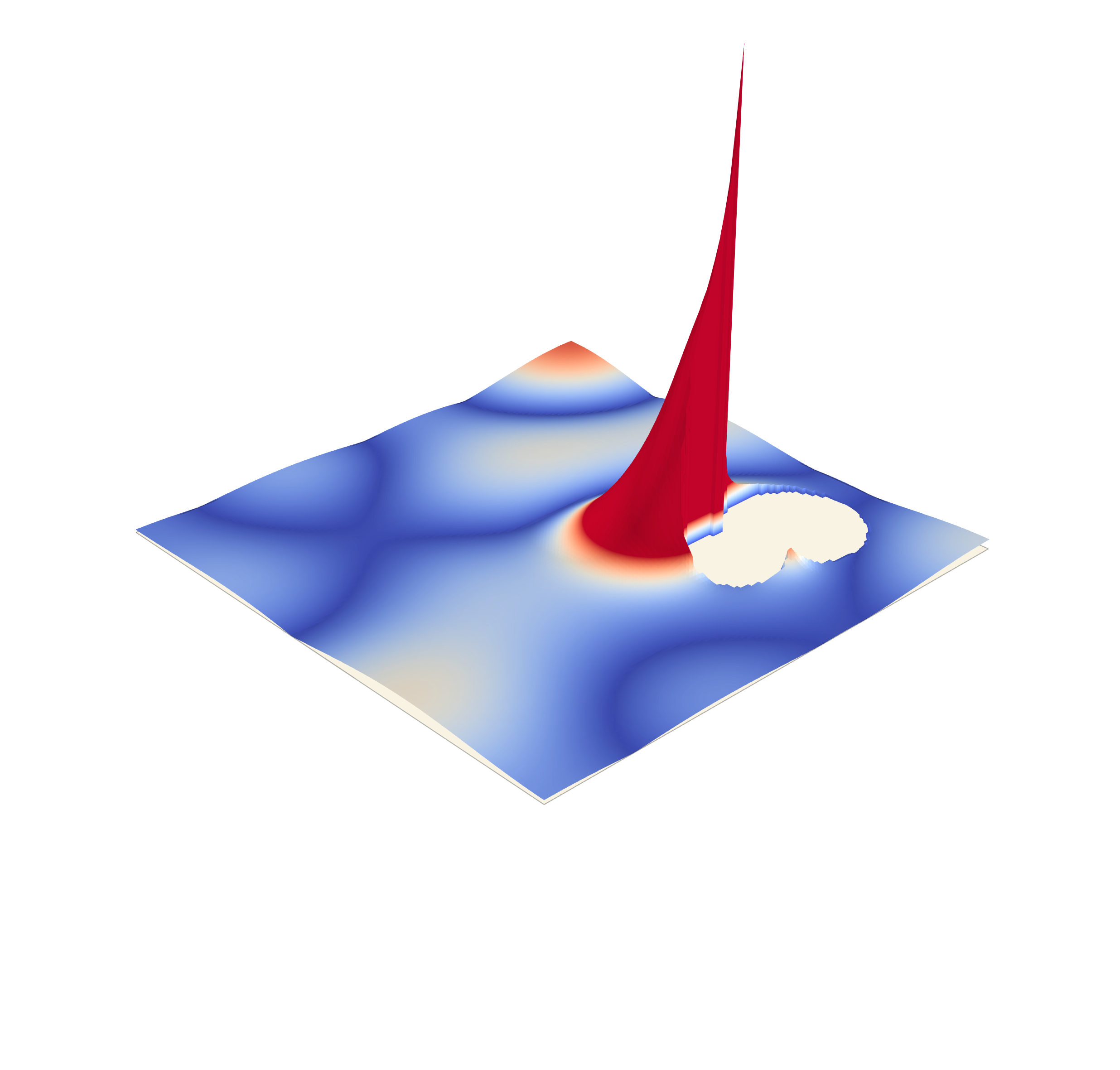}}} 
\subfigure[Using proposed method]{{\includegraphics[trim={10cm 17cm 10cm 6cm}, clip, width=.25\textwidth]{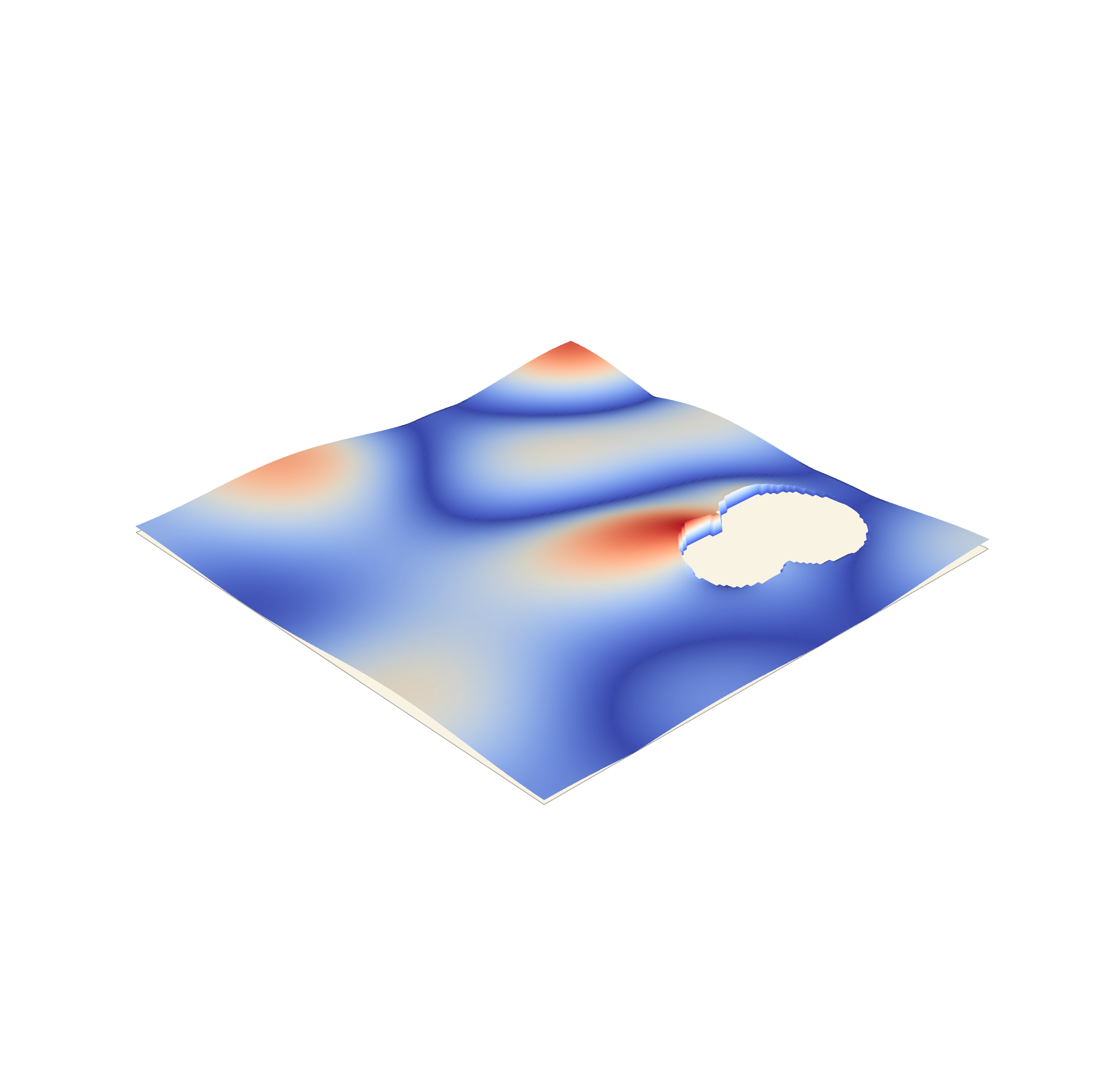}}} 
\end{center}
\vspace{-.3cm}\vspace{-.3cm}\caption{\it Error distribution at the final time moment ($t=1$) in case of a non-smooth moving object using different extrapolation approaches.} \label{fig::diffusion_union_error}
\end{figure}

\begin{figure}[!h]
\begin{center}
\includegraphics[width=.9\textwidth]{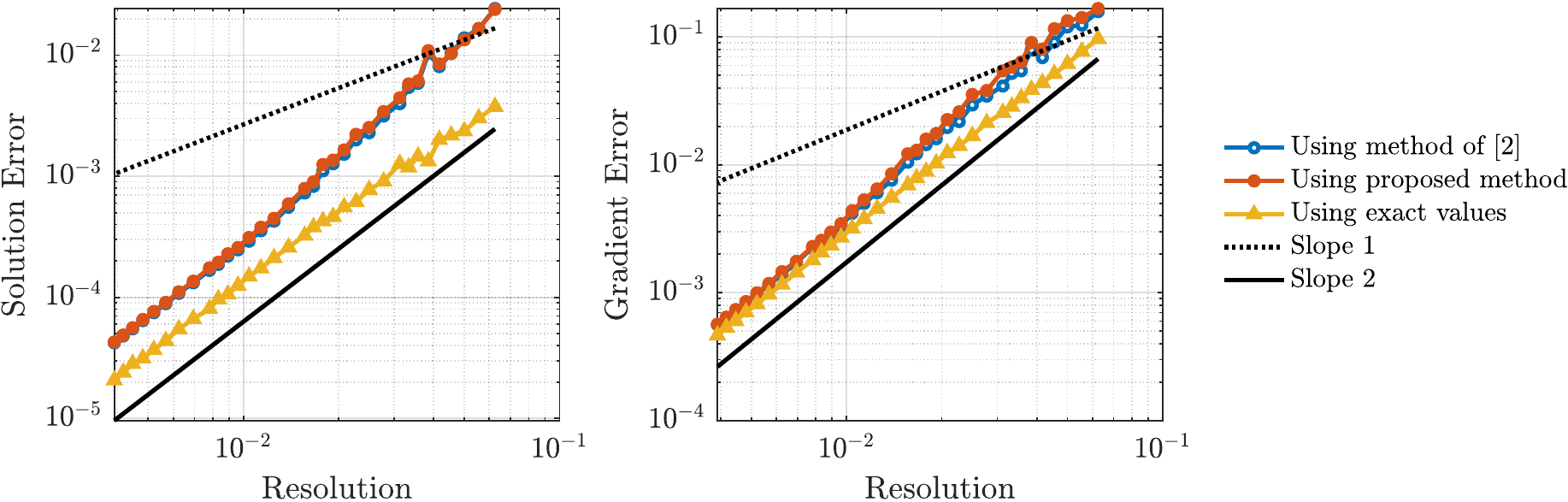}
\end{center}
\vspace{-.3cm}\vspace{-.3cm}\caption{\it Accuracy of solving diffusion equation (in the $L^\infty$ norm) in case of a smooth moving object using different extrapolation approaches.} \label{fig::diffusion_circle_convergence}
\end{figure}

\begin{figure}[!h]
\begin{center}
\includegraphics[width=.9\textwidth]{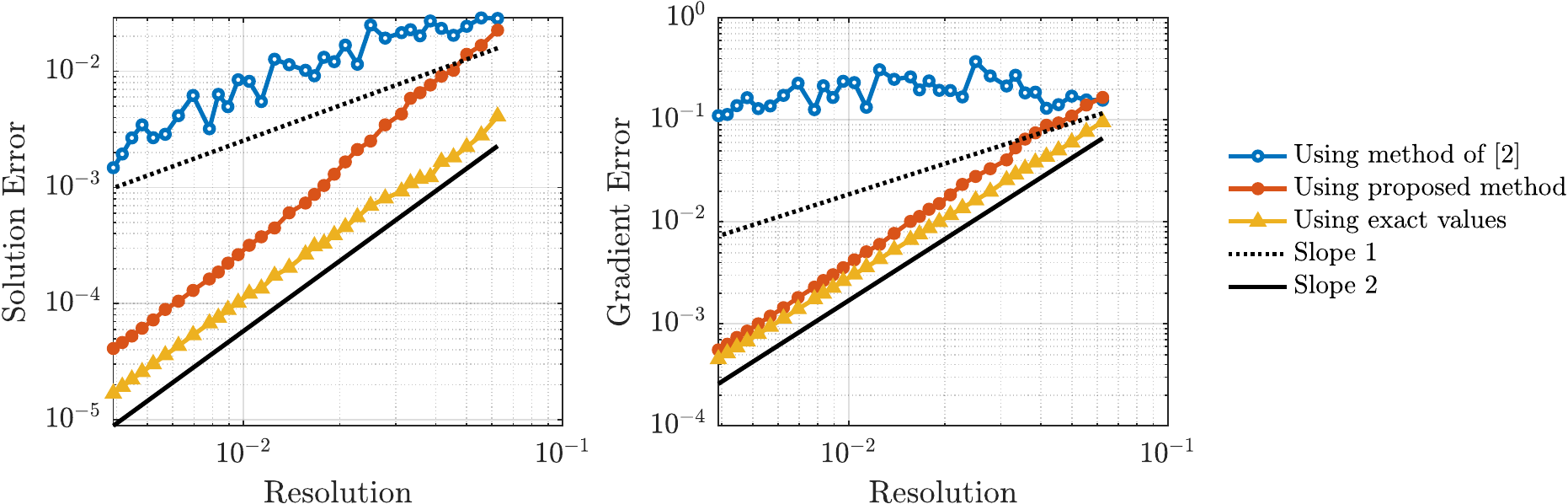}
\end{center}
\vspace{-.3cm}\vspace{-.3cm}\caption{\it Accuracy of solving diffusion equation (in the $L^\infty$ norm) in case of a non-smooth moving object using different extrapolation approaches.} \label{fig::diffusion_union_convergence}
\end{figure}

\reviewerOne{
Obtained results are summarized in Figures \ref{fig::diffusion_circle_error}, \ref{fig::diffusion_union_error}, \ref{fig::diffusion_circle_convergence} and \ref{fig::diffusion_union_convergence}. Figures \ref{fig::diffusion_circle_error} and \ref{fig::diffusion_union_error} show error distributions at the final time $t=1$ on a $128^2$ grid. As one can see, for a smooth moving object (Fig. \ref{fig::diffusion_circle_error}) using either the ND-PDE extension or the WCD-PDE one, results in errors that are very close to ones of the reference solution, while for a non-smooth moving object (Fig. \ref{fig::diffusion_union_error}) the WCD-PDE extension produces a much more accurate solution that is also very close to the reference one, whereas using the ND-PDE extension results in a significant accumulation of errors in the wake of the moving object. More quantitative conclusions can be drawn from the convergence studies shown in Figures \ref{fig::diffusion_circle_convergence} and \ref{fig::diffusion_union_convergence}, in which the dependence of the error in the $L^\infty$ norm for both numerical solution and its gradient on the grid resolution is presented. In case of a smooth moving object (Fig. \ref{fig::diffusion_circle_convergence}) both extension methods lead to second-order convergence both in the numerical solutions and its gradient (as what is expected from the superconvergent method of \cite{Bochkov;Gibou:19:Solving-the-Poisson-}) with magnitude of errors being very close to the ones of the reference solution. In case of a non-smooth moving object (Fig. \ref{fig::diffusion_union_convergence}) the accuracy of numerical solutions obtained using the ND-PDE extension degrades severely showing only first-order convergence in the solution itself and non-convergence in the gradient. At the same time the accuracy of computations based on the proposed WCD-PDE extension seems affected very slightly by the presence of sharp corners, retaining the second-order convergence in the numerical solution and its gradient with errors very close to the reference solution.
}

\reviewerOne{
\textbf{Remark:} Note that the presented in this work extrapolation approach is designed for extending smooth scalar fields (as in the present example). However, in general, solutions to partial differential equations in domains with sharp features may contain singularities. In such cases for best results the proposed extrapolation procedure should only be applied to the regular part of the solution, the singular part must be dealt with separately using special methods, for example, as in \cite{Wigley:88:An-efficient-method-}.
}

\section{Conclusion} \label{sec:Conclusion}
We have presented a numerical method for extrapolating scalar quantities across the boundaries of irregular domains that may present high-curvature features or kinks. Linear and quadratic extrapolations procedures produce second- and third-order accurate results in the $L^\infty$ norm, respectively and do so regardless of the irregularity of the boundaries, i.e. boundaries with kinks can readily be considered. These procedures are effective in both two and three spatial dimensions and can be implemented on quadtree and octree Cartesian grids. We have shown through numerical examples that errors associated with extrapolations can be reduced by \mylinelabel{rev1:typo}\reviewerOne{several orders of magnitude} in some cases, compared with the \reviewerOne{approach of \cite{Aslam:04:A-partial-differenti} commonly} used in level-set methods. \mylinelabel{rev1:example_outro}\reviewerOne{We have also presented an example of solving a diffusion equation on evolving domains in order to highlight the importance of accurate extrapolation near sharp geometric features for practical applications.} The numerical method we introduced is based on solving PDEs in pseudo-time, but we note that static solutions based on an implicit approach like Fast Marching or Fast Sweeping could be obtained and we expect the results to follow the same general behavior as that presented in the current manuscript.

\section*{Acknowledgment} This work was supported by ONR MURI N00014-17-1-2676.
\newpage
\clearpage
\setcounter{page}{1}
\bibliographystyle{plain}
\addcontentsline{toc}{section}{\refname}
\bibliography{references_short}

\end{document}